\documentclass[12pt]{article}

\usepackage[latin1]{inputenc}
\usepackage{bbm}
\usepackage{stmaryrd}
\usepackage{latexsym}
\usepackage{amsfonts}
\usepackage{amsmath,amssymb,amscd}
\usepackage{yfonts}
\usepackage{dsfont}
\usepackage{mathabx}
\usepackage[dvips]{graphicx}
\usepackage{epsfig}
\usepackage{psfrag}
\usepackage[font={small}]{caption}

\DeclareFontFamily{U}{txsyc}{}
\DeclareFontShape{U}{txsyc}{m}{n}{
   <-> txsyc%
}{}
\DeclareFontShape{U}{txsyc}{bx}{n}{
   <-> txbsyc%
}{}
\DeclareFontShape{U}{txsyc}{l}{n}{<->ssub * txsyc/m/n}{}
\DeclareFontShape{U}{txsyc}{b}{n}{<->ssub * txsyc/bx/n}{}
\DeclareSymbolFont{symbolsC}{U}{txsyc}{m}{n}
\SetSymbolFont{symbolsC}{bold}{U}{txsyc}{bx}{n}
\DeclareFontSubstitution{U}{txsyc}{m}{n}
\DeclareMathSymbol{\df}{\mathrel}{symbolsC}{"42}
\DeclareMathSymbol{\fd}{\mathrel}{symbolsC}{"43}
\DeclareMathSymbol{\lJoin}{\mathrel}{symbolsC}{"58}
\DeclareMathSymbol{\rJoin}{\mathrel}{symbolsC}{"59}

\newcommand{\cC}{{\cal C}}

\newcommand{\cE}{{\cal E}}
\newcommand{\cF}{{\cal F}}

\newcommand{\cL}{{\cal L}}

\newcommand{\cO}{{\cal O}}
\newcommand{\cP}{{\cal P}}

\newcommand{\CC}{\mathbb{C}}

\newcommand{\EE}{\mathbb{E}}

\newcommand{\LL}{\mathbb{L}}
\newcommand{\NN}{\mathbb{N}}
\newcommand{\PP}{\mathbb{P}}

\newcommand{\RR}{\mathbb{R}}

\newcommand{\ZZ}{\mathbb{Z}}

\newcommand{\iy}{\infty}
\newcommand{\lt}{\left}
\newcommand{\me}{\medskip}

\newcommand{\pa}{\partial}
\newcommand{\ri}{\rightarrow}
\newcommand{\rt}{\right}
\newcommand{\sm}{\smallskip}

\newcommand{\wi}{\widetilde}
\newcommand{\wit}{\widehat}

\newcommand{\essinf}{\mathrm{essinf}}
\newcommand{\esssup}{\mathrm{esssup}}
\newcommand{\ex}{\exists\ }
\newcommand{\fo}{\forall\ }

\newcommand{\lan}{\lt\langle}
\newcommand{\lve}{\lt\vert}
\newcommand{\lVe}{\lt\Vert}

\newcommand{\ran}{\rt\rangle}
\newcommand{\rve}{\rt\vert}
\newcommand{\rVe}{\rt\Vert}

\newcommand{\st}{\,:\,}

\newcommand{\un}{\mathds{1}}

\newcommand{\bq}{\begin{eqnarray*}}
\newcommand{\bqn}[1]{\begin{eqnarray}\label{#1}}
\newcommand{\eq}{\end{eqnarray*}}
\newcommand{\eqn}{\end{eqnarray}}
\newcommand{\wwtbp}{\par\hfill $\blacksquare$\par\me\noindent}
\newcommand{\thistitlepagestyle}{}
\newcommand{\lin}{\llbracket}
\newcommand{\rin}{\rrbracket}

\newcommand{\ttsim}{\raise.17ex\hbox{$\scriptstyle\mathtt{\sim}$}}

\newtheorem{pro}{Proposition} 

\newtheorem{lem}[pro]{Lemma}
\newtheorem{theo}[pro]{Theorem}
\renewcommand{\thepro}{\arabic{pro}}

\newenvironment{rem}
{\par\me\refstepcounter{pro}\noindent{\bf Remark \thepro\ }}
{\par\hfill $\square$\par\me\noindent}
\newenvironment{rems}
{\par\me\refstepcounter{pro}\noindent{\bf Remarks \thepro\ }\par\sm}
{\par\hfill $\square$\par\me\noindent}
\newcommand{\proof}{\par\me\noindent\textbf{Proof}\par\sm\noindent}

\title{On quantitative convergence to quasi-stationarity}
\author{Persi Diaconis${}^\dagger$ and Laurent Miclo${}^\ddagger$
}
\date{\box1
 \box2
}

\begin{document}

\setbox1=\vbox{
\large
\begin{center}
${}^\dagger$Department of Mathematics\\
Stanford University, California, USA\\
\end{center}
} 
\setbox2=\vbox{
\large
\begin{center}
${}^\ddagger$Institut de Mathématiques de Toulouse, UMR 5219\\
Universit\'e de Toulouse and CNRS, France\\
\end{center}
} 
\setbox3=\vbox{
\hbox{Stanford University\\}
\hbox{Department of Mathematics\\}
\hbox{Building 380, Sloan Hall\\}
\hbox{Stanford, California 94305, USA\\}
}
\setbox4=\vbox{
\hbox{${}^\ddagger$miclo@math.univ-toulouse.fr\\}
\vskip1mm
\hbox{Institut de Mathématiques de Toulouse\\}
\hbox{Université Paul Sabatier\\}
\hbox{118, route de Narbonne\\} 
\hbox{31062 Toulouse Cedex 9, France\\}
}
\setbox5=\vbox{
\box3
\vskip5mm
\box4
}

\maketitle
\thistitlepagestyle
\abstract{
The quantitative long time behavior of absorbing, finite, irreducible Markov processes is considered.
Via Doob transforms, it is shown that only the knowledge of the ratio of
the values of the underlying first Dirichlet eigenvector is necessary to come back to the well-investigated situation
of the convergence to equilibrium of ergodic finite Markov processes.
This leads to explicit estimates on the convergence to quasi-stationarity, in particular via
functional inequalities. When the process is reversible, the optimal exponential rate consisting of the spectral gap between 
the two first Dirichlet eigenvalues is recovered.
Several simple examples are provided to illustrate the  bounds obtained.
}
\vfill\null
{\small
\textbf{Keywords: }
Absorbing finite Markov process, quasi-stationary measure, Doob transform, quantitative bounds on convergence,
spectral gap, logarithmic Sobolev inequality, birth and death process.
\par
\vskip.3cm
\textbf{MSC2010:} primary: 60J27, secondary: 60E15, 46E39, 47D08, 37A30, 15B51.
}\par

\newpage 

\section{Introduction}\label{intro}

This paper begins to develop a quantitative theory of rates of convergence to quasi-stationarity, as in the following example.
Consider the simple symmetric random walk on $\{0, ..., N\}$ with holding $1/2$ at $N$ and absorbing at 0. Let $\bar X_t$
 be the position of the walk at time $t\in\ZZ_+$ and $T$ be the absorption time at 0. Let $\mu_t(x)\df \PP[\bar X_t=x\vert T>t]$, for $x\in\{1, ..., N\}$.
 Classical theory, reviewed below, shows that 
 \bq
 \lim_{t\ri +\iy}\mu_t(x)&=& \nu(x)\ \df\ Z^{-1}\cos\lt(\frac{(2N+1-2x)\pi}{2(2N-1)}\rt)\eq
 with $Z^{-1}\df
  2\tan\lt(\frac{\pi}{2(2N-1)}\rt)$, the normalizing constant. The measure $\nu$ is called a quasi-stationary distribution. How large does $t$ have to be so that these asymptotics are useful? 
  In Section \ref{examples}, which is devoted to explicit computations, we prove for the continuous time counterpart of the above process that
  for any starting distribution on $\{1, ..., N\}$ and for all $s\geq 0$,
  \bqn{exa0}
  \lVe \mu_t-\nu\rVe_{\mathrm{tv}}&\leq &\frac{2\sqrt{2}}{\pi^2}(1+\cO(N^{-1}))\exp(-s)\eqn
  \bq
t&=&\frac5{2\pi^2}N^2\ln(N)+\frac{s}{\pi^2}N^2\eq
Thus the quasi-stationary asymptotics takes hold for $t$ larger than $N^2\ln(N)$.
In \eqref{exa0}, $\mu_t$ and $\nu$ depend on $N$ but the bounds are uniform in $N$.\par\me
We will work mainly in the continuous time setting, which is more convenient to deal with. We will come back to the discrete time framework in Section \ref{sdtm}.
Generally, a   quasi-stationary distribution of an absorbing Markov process $\bar X\df(\bar X_t)_{t\geq 0}$ is a probability measure $\nu$ on the  state space $S$ (where the absorbing points have been removed)
such that starting from this distribution, the time marginal laws $\cL(\bar X_t)$  remain proportional to $\nu$ on $S$, for all $t\geq 0$.
For nice processes $\bar X$, the quasi-stationary distribution is unique and starting from any distribution on $S$, the conditional (to non-absorption) law
$\mu_t\df\cL(\bar X_t\vert \bar X_t\in S)$ converges toward $\nu$ for large times $t\geq 0$.
The purpose of this article is to investigate this convergence  quantitatively when
$S$ is finite.\par\me
More precisely, the framework is as follows. The whole finite state space is $\bar S\df S\sqcup\{\iy\}$, where $\iy$ is the absorbing point.
There is no loss of generality in assuming there is only one such point, up to lumping together all the absorbing points.
Let $\bar L$ be the generator of the process $\bar X$ on $\bar S$, seen as a matrix $(\bar L(x,y))_{x,y\in\bar S}$. 
To any given probability measure $m_0$ on $\bar S$, there is a unique (in law) Markov process $\bar X$  whose generator is $\bar L$ and whose 
initial distribution $\cL(\bar X_0)$ is $m_0$.
For any $t\geq 0$, let $m_t=\cL(\bar X_t)$. Using matrix notation, where measures are seen as row vectors (and functions as column vectors),
we have
\bq
\fo t\geq 0,\qquad m_t&=&m_0\bar P_t\eq
where $(\bar P_t)_{t\geq 0}$ is the semi-group $(\exp(t\bar L))_{t\geq 0}$ associated to $\bar L$.
Except if $m_0$ is the Dirac mass on $\iy$, for any $t\geq 0$, $m_t(S)>0$ and we can define the probability measures $\mu_t$ as
the restriction to $S$ of $m_t/m_t(S)$. 
They will be our main objects of interest here.
By definition, we have
\bqn{mut}
\fo t\geq 0,\,\fo f\in \cF,\qquad \mu_t[f]&=&\frac{\mu_0 [\bar P_t[f]]}{\mu_0[ \bar P_t[\un_S]]}\eqn
where  $\cF$ is the space of real functions defined on $S$, also seen as functions defined on $\bar S$ 
which vanish at $\iy$ 
(Dirichlet condition at $\iy$).
A probability measure $\nu$ on $S$ is said to be a quasi-stationary measure for $\bar L$ if $\mu_0=\nu$ implies that $\mu_t=\nu$
for all $t\geq 0$. We will recall below a convenient assumption ensuring there is a unique quasi-invariant measure $\nu$
associated to $\bar L$. The objective 
of this paper is to quantify the convergence of $\mu_t$ toward $\nu$ for large times $t\geq 0$, whatever the initial distribution $\mu_0$.\par\sm
For any 
$x\in S$, denote $V(x)=L(x,\iy)\geq 0$, the killing rate at $x$. The symbol $V$ will designate the function $S\ni x\mapsto V(x)$
as well as the $S\times S$ diagonal matrix whose values on the diagonal are given by $V$, namely the multiplication operator by $V$
on $\cF$.
Let $L$ be the Markov generator on $S$ which is such that the $S\times S$ minor of $\bar L$ can be written $L-V$.
Our main assumption is that $L$ is irreducible. At some point, this hypothesis will 
be	
strengthened by a reversibility assumption, in order to get more explicit results.
A traditional application of the Perron-Frobenius theorem 
(see for instance the book \cite{MR2986807} of Collet,  Mart{\'{\i}}nez and San
              Mart{\'{\i}}n)
to $L-V$ or to the associated semi-group, seen as operators on measures on $S$,
ensures that there exists a unique quasi-invariant measure $\nu$ associated to $\bar L$. The probability measure $\nu$ gives a positive weight to any point of $S$.
Furthermore there exists $\lambda_1\geq 0$ such that $\nu (L-V)=-\lambda_1\nu$, $\lambda_1$ is the eigenvalue of $V-L$ which is strictly less than  the real parts of the remaining eigenvalues (in $\CC$).
In the same manner, there exists a unique invariant measure $\eta$ for $L$, charging all points of $S$. To see the relation between $\nu$ and $\eta$,
consider the operator $L^*$ which is adjoint to $L$ in $\LL^2(\eta)$.
As a matrix, it is given by
\bq
\fo x,y\in S,\qquad L^*(x,y)&=&\frac{\eta(y)}{\eta(x)}L(y,x)\eq
The fact that $\eta$ is invariant is equivalent to the fact that $L^*$ is a Markovian generator.
We can thus
apply the Perron-Frobenius theorem to $L^*-V$, seen as an operator on $\cF$ to find a positive function $\varphi^*$ on $S$ such that
$(L^*-V)[\varphi^*]=-\lambda_1\varphi^*$. Let us renormalized $\varphi^*$ so that $\eta[\varphi^*]=1$.
Then $\nu=\varphi^*\cdot\eta$, the probability measure admitting the density $\varphi^*$ with respect to $\eta$.
Indeed, for any test function $f\in\cF$, we have
\bq
(\varphi^*\cdot\eta)[(L-V)[f]]&=&\eta[\varphi^*(L-V)[f]]\\
&=&\eta[(L^*-V)[\varphi^*]f]\\
&=&-\lambda_1\eta[\varphi^* f]\\
&=&-\lambda_1(\varphi^*\cdot\eta)[f]\eq
so that $(\varphi^*\cdot\eta)(L-V)=-\lambda_1(\varphi^*\cdot\eta)$ and by consequence $(\varphi^*\cdot\eta)\bar P_t=\exp(-\lambda_1 t)(\varphi^*\cdot\eta)+(1-\exp(-\lambda_1 t))\delta_\iy$.\par
This relation implies that if the process $\bar X$ is started from the quasi-distribution $\nu$, then the absorption time $\tau\df\inf\{t\geq 0\st \bar X_t=\iy\}$ is distributed as an exponential distribution
of parameter $\lambda_1$. Indeed, we have for any $t\geq 0$,
\bq
\PP_\nu[\tau>t]&=&\nu\bar P_t[S]\\
&=&\exp(-\lambda_1t)\eq
where $\PP_{\nu}$ is the underlying probability measure, when $\bar X_0$ is distributed according to $\nu$.
More generally, from this identity, it is not difficult to deduce that for any initial distribution $m_0$ not equal to $\delta_\iy$, we have
\bq
\lim_{t\ri+\iy}\ln(\PP_{m_0}[\tau >t])&=&-\lambda_0\eq
showing that $\lambda_0$ is the exponential rate of absorption.
\par
Furthermore, we can find a positive function $\varphi\in \cF$ such that 
$(L-V)\varphi=-\lambda_1\varphi$, but we rather normalize it through the relation $\eta[\varphi^2]=1$.
For any positive function $f\in\cF$, we note $f_\wedge\df\min_{x\in S}f(x)$ and $f_\vee\df\max_{x\in S}f(x)$.\par
Finally, consider the Markovian operator $\wi L$ on $S$ which is defined by its off-diagonal entries via
\bqn{wiL}
\fo x\not=y\in S,\qquad \wi L(x,y)&\df& L(x,y)\frac{\varphi(y)}{\varphi(x)}\eqn
Let $(\wi P_t)_{t\geq 0}$ be the associated Markovian semi-group.
Since $\wi L$ is irreducible, it admits an invariant probability $\wi\eta$. In next section we will check that it is given by
\bqn{wieta}
\fo x\in S,\qquad \wi \eta(x)&=&\frac{\varphi(x)\varphi^*(x)\eta(x)}{\sum_{y\in S}\varphi(y)\varphi^*(y)\eta(y)}\eqn
\par
To give a first estimate on the convergence of $\mu_t$ toward $\nu$, let us recall that the total variation of a signed measure $m$ on $S$ satisfying $m(S)=0$
is given equivalently by
\bq
\lVe m\rVe_{\mathrm{tv}}&\df& 2\sup_{A\subset S}\lve m(A)\rve\\
&=&\sup_{f\in\cF,\, \lVe f\rVe_\iy\leq 1}m(f)\\
&=&\sum_{x\in S}\lve m(x)\rve\eq
(where as usual, $\lVe f\rVe_\iy$ designates the supremum norm of $f$). Note this definition differs by a factor of 2 from the probabilist version.
\par
\begin{theo}\label{tvtv}
For any probability measure $\mu_0$ on $S$ and for any $t\geq 0$, we have
\bq
 \frac{\varphi_\wedge}{2\varphi_\vee}\lVe \wi\mu_0\wi P_t -\wi\eta\rVe_{\mathrm{tv}}\ \leq \ \lVe \mu_t-\nu\rVe_{\mathrm{tv}}\ \leq \ 2\frac{\varphi_\vee}{\varphi_\wedge}\lVe \wi\mu_0\wi P_t -\wi\eta\rVe_{\mathrm{tv}}\eq
where $\wi\mu_0$ is the probability on $S$ whose density with respect to $\mu_0$ is proportional to $\varphi$.
In particular the asymptotic exponential rate of convergence of $\lVe \mu_t-\nu\rVe_{\mathrm{tv}}$ and 
$\lVe \wi\mu_0\wi P_t -\wi\eta\rVe_{\mathrm{tv}}$ are the same.
\end{theo}
\par
Note that in the trivial case where there is no absorption, namely $V\equiv 0$, we have $\varphi\equiv 1\equiv\varphi^*$, $(\wi P_t)_{t\geq 0}=( P_t)_{t\geq 0}$,
the Markovian semi-group generated by $L$,  $\nu=\wi\eta$ and $\mu_t=\mu_0\wi P_t$ for all $t\geq 0$, so that the above bounds are optimal, up to the factor 2.
\par
In a forthcoming paper, we investigate the quantity $\varphi_\vee/\varphi_\wedge$, providing different upper bounds  via path and spectral considerations,
first step toward the extension of the results presented here to certain denumerable chains.
\par\sm
Theorem \ref{tvtv} reduces the study of convergence to quasi-stationarity to the much more well-studied situation of the convergence to equilibrium.
One can for instance resort to functional inequality techniques (see for instance the lecture notes of Saloff-Coste \cite{MR1490046}), the simplest of them being the  $\LL^2$ approach. Let $\wi L^\diamond$ be the additive symmetrization of $\wi L$ in $\LL^2(\wi\eta)$: it is equal to $(\wi L+\wi L^*)/2$, where
$\wi L^*$ is the adjoint operator of $\wi L$ in $\LL^2(\wi\eta)$. The matrix of this Markov generator is described by its off-diagonal entries:
\bqn{Ldiamond}
\fo x\not=y\in S, \qquad
\wi L^\diamond(x,y)&=&\frac12\lt(L(y,x)\frac{\varphi^*(y)\eta(y)}{\varphi^*(x)\eta(x)}+L(x,y)\frac{\varphi(y)}{\varphi(x)}\rt)\eqn
Its self-adjointness implies that $\wi L^\diamond$ is diagonalizable in $\RR$ and let $\wi\lambda> 0$ stands for
the smallest non-zero eigenvalue of $-\wi L^\diamond$. Since $\wi L^\diamond$ is irreducible, the eigenvalue 0 has multiplicity 1 (with
eigenspace consisting of the constant functions) and $\wi\lambda$ is the spectral gap of $\wi L^\diamond$.
Then we get:
\begin{theo}\label{tvL2}
For 
any $t\geq 0$, we have
\bq
\sup_{\mu_0\in\cP}\lVe \mu_t-\nu\rVe_{\mathrm{tv}}&\leq &
\sqrt{\frac{\eta[\varphi\varphi^*]}{(\varphi\varphi^*\eta)_{\wedge}}}
\frac{\varphi_\vee}{\varphi_\wedge}\exp(-\wi\lambda t)\eq
where $\cP$ stands for the set of probability measures on $S$.
\end{theo}
\par
We have $\wi\lambda=A^{-1}$, where $A$ is the smallest positive constant such that the following Poincaré inequality is satisfied for all $f\in\cF$,
\bqn{Poincare}
\sum_{x\in S}(f(x)-\wi\eta[f])^2\, \varphi^*(x)\varphi(x)\eta(x)&\leq &A \sum_{x,y\in S}(f(y)-f(x))^2\, \varphi^*(x)\varphi(y)\eta(x) L(x,y)\eqn
This variational formulation enables to compare $\wi\lambda$ with $\lambda$ (see for instance Diaconis and Saloff-Coste \cite{MR1233621} and Fill \cite{MR1097464}),
the spectral gap of the additive symmetrization of $L$ in $\LL^2(\mu)$:
\bqn{will}
\wi\lambda&\geq & \frac{\varphi_\wedge\varphi^*_\wedge}{\varphi_\vee\varphi^*_\vee}\lambda\eqn
 We will put these considerations  into practice in  Example \ref{anre}.
\par\sm
Let us now assume that $\eta$ is reversible for $L$. Then $-(L-V)$ is self-adjoint in $\LL^2(\eta)$
and so is diagonalizable in $\RR$. As it was already mentioned for the general case, its smallest eigenvalue is $\lambda_1>0$.
Consider its next eigenvalue $\lambda_2>\lambda_1$ (the strict inequality is a consequence of the irreducibility of $L$ in the Perron-Frobenius theorem).
The next result shows that to get a useful understanding of the convergence of $\mu_t$ toward $\nu$ for large $t\geq 0$,
only the knowledge of $\eta$, of the ratio of the extrema of $\varphi$ and of $\lambda_2-\lambda_1$ is required.
\begin{theo}\label{tvL2rev}
Under the reversibility assumption, 
 for any $t\geq 0$, we have
\bq
\sup_{\mu_0\in\cP}\lVe \mu_t-\nu\rVe_{\mathrm{tv}}&\leq &
\sqrt{\frac{1}{(\varphi^2\eta)_{\wedge}}}
\frac{\varphi_\vee}{\varphi_\wedge}\exp(-(\lambda_2-\lambda_1) t)\\
&\leq & \sqrt{\frac{1}{\eta_{\wedge}}}
\lt(\frac{\varphi_\vee}{\varphi_\wedge}\rt)^2\exp(-(\lambda_2-\lambda_1) t)
\eq
\end{theo}
\par
Note that \eqref{mut} can be written in terms of Feynman-Kac integrals. Let $(X_t)_{t\geq 0}$ be a Markov process starting from the initial law $\mu_0$
and admitting $L$ as generator. We have
\bq
\fo t\geq 0,\,\fo f\in \bar\cF,\qquad \mu_t[f]&=&\frac{\EE_{\mu_0}\lt[f(X_t)\exp\lt(-\int_0^t V(X_s)\, ds\rt)\rt]}{\EE_{\mu_0}\lt[\exp\lt(-\int_0^t V(X_s)\, ds\rt)\rt]}\eq
The stability for large times of such expressions have been extensively studied by Del Moral and his coauthors (see for instance his recent book \cite{MR3060209} and the references given there).
They also use  estimates on the convergence to equilibrium of Markov processes. Since their  assumptions are based on Dobrushin type conditions
on the underlying Markov process (or on some of its modifications, see e.g.\ Del Moral and Miclo \cite{MR1956078}), the deduced bounds are often quite coarse.
While we work in the same spirit, we will rather resort to spectral techniques, which lead to relatively sharp estimates, as will be illustrated by several examples.
In particular, we obtain in the reversible case the optimal asymptotical rate $\lambda_2-\lambda_1$ (see e.g.\ the review paper of M{\'e}l{\'e}ard and Villemonais \cite{MR2994898}, with a non-quantified pre-exponential factor).
Under appropriate conditions, this rate was deduced asymptotically  for birth and death processes by van Doorn
\cite{MR1133722} (see also van Doorn and Zeifman
\cite{MR2576013} for another example), which are outside the scope of the present note, because the state space is not finite.
We hope that in a future work, we will be able to extend the above quantitative bounds to  more general situations 
of appropriate denumerable Markov processes or diffusions, perhaps under the condition there is a unique quasi-invariant measure (usually this requires that the process 
comes in from infinity  fast enough, see for instance  Collet,  Mart{\'{\i}}nez and San
              Mart{\'{\i}}n \cite{MR2986807}). For Brownian motion absorbed on the boundary of a compact domain in Euclidean spaces, one may see
      Gyrya       and Saloff-Coste \cite{MR2807275} and
              Lierl and Saloff-Coste \cite{2012arXiv1210.4586L}
              \par
  \sm
  The literature on quasi-stationarity is substantial and we are able to call on several comprehensive surveys. One short readable survey, close in spirit to our paper, is by Van Doorn and Pollett \cite{MR3063313} (discrete state space, continuous time). More general state spaces and applications in biology are emphasized by Méléard and Villemonais \cite{MR2994898}. A recent book length treatment by 
  Collet,  Mart{\'{\i}}nez and San
              Mart{\'{\i}}n  \cite{MR2986807} treats all aspects. All of these review the history (Yaglom, Bartlett, Darroch-Seneta, ...). A most useful adjunct to these surveys is the annotated online bibliography kept up to date by Phil Pollett, see \texttt{http://www.maths.uf.edu.au/$\mathbf{\sim}$pkp/papers/qsds.html}.
\par
We have not found very much literature on the kind of quantitative questions treated here. A useful review of previous quantitative efforts is in Section 4 of Van Doorn and Pollett \cite{MR3063313}. This is along the lines of spectral gap estimates without consideration of the size of the state space or  the starting distribution.
Some quantitative bounds are also deduced in  the recent papers of Cloez and Thai \cite{2013arXiv1312.2444C} and of Champagnat and Villemonais \cite{2014arXiv1404.1349C}.            \par\me
The plan of the paper is very simple: the next section presents the proof of the above theorems, as well as an alternative bound based on logarithmic Sobolev inequalities,
 the next section contains some illustrative examples. The final section gives further examples in discrete time.

\section{Proofs}

The following arguments are based on a simple use of Doob's transforms, which by a conjugation by $\varphi$,
replace $V$ by a constant killing rate.

\subsection{Proof of Theorem \ref{tvtv}}

Let $\Phi$ be the diagonal matrix corresponding to the multiplication by $\varphi$ operating on $\cF$.
Thus $\Phi^{-1}$ is just the diagonal matrix corresponding to the multiplication by $1/\varphi$.
We begin by checking that the generator matrix $\wi L$ defined in (\ref{wiL})
satisfies
\bqn{wiLL}
\wi L&=&\Phi^{-1}(L-V+\lambda_1I)\Phi\eqn
where $I$ is the identity matrix.
 Indeed, the off-diagonal entries of the r.h.s.\ coincide with those of $\Phi^{-1}L\Phi$
 which are those of $\wi L$ by (\ref{wiL}).
Thus it is sufficient to check that the sums of the rows of $\Phi^{-1}(L-V+\lambda_1I)\Phi$ vanish.
The sum corresponding to the row indexed by $x\in S$ is
\bq
\frac1{\varphi(x)}(L[\varphi](x)-V(x)\varphi(x))+\lambda_1
&=&0\eq
since by definition, $\varphi$ is an eigenfunction of $L-V$ associated to the eigenvalue $-\lambda_1$.
\par
It is now easy to check \eqref{wieta}: it must be seen that
\bq
\fo f\in\cF,\qquad \wi \eta[\wi L[f]]&=&0\eq
From \eqref{wiLL}, the l.h.s.\ is equal
to
\bq
\wi \eta[\varphi^{-1}(L-V+\lambda_1)[\varphi f]]&=&\eta[\varphi^*(L-V+\lambda_1)[\varphi f]]/\eta[\varphi\varphi^*]\\
&=&\eta[\varphi f(L^*-V+\lambda_1)[\varphi^*]]/\eta[\varphi\varphi^*]\\
&=&0\eq
because $\varphi^*$ is an eigenfunction of $L^*-V$ associated to the eigenvalue $-\lambda_1$.\par
Next rewrite \eqref{wiLL} in the form 
\bqn{conj}
\Phi (\wi L-\lambda_1I)\Phi^{-1}&=&L-V\eqn 
and exponentiate this identity to find
\bq
\fo t\geq 0,\qquad \exp(-\lambda_1 t)\Phi \wi P_t\Phi^{-1}&=&\bar P_t\eq
(the r.h.s.\ is to be understood as the restriction of $\bar P_t$ to $\cF$, as explained after \eqref{mut}).
Thus for any $\mu_0\in\cP$ and $f\in\cF$, we have
\bq
\fo t\geq 0,\qquad \exp(-\lambda_1 t)\mu_0[\varphi]\wi\mu_0[\wi P_t[f/\varphi]]&=&\mu_0[\bar P_t[f]]\eq
(recall that $\wi\mu_0$ is the probability on $S$ whose density with respect to $\mu_0$ is proportional to $\varphi$).
We deduce from \eqref{mut} that
\bqn{mutwimu0}
\fo t\geq 0,\qquad \mu_t[f]&=&\frac{\wi\mu_0[\wi P_t[f/\varphi]]}{\wi\mu_0[\wi P_t[1/\varphi]]}\eqn
Since $\wi P_t$ converges to $\wi\eta$ as $t$ goes to infinity, we get  that
\bq
\lim_{t\ri+\iy} \mu_t[f]&=&\frac{\wi\eta[f/\varphi]}{\wi \eta[1/\varphi]}\\
&=&\nu[f]\eq
due to the proportionality between the measures $\nu$, $\varphi^*\cdot \eta$ and $\varphi^{-1}\cdot\wi\eta$.
Thus the convergence toward quasi-stationarity has been recovered.\par
To get an estimate on the speed of convergence, 
we need the two following basic lemmas. \par
On a general measurable space, consider two probability measures $\wi\mu\ll\wi\nu$, as well as a measurable function
$\psi>0$. Define
\bq
\mu\,\df\, \frac{\psi}{Z_{\wi\mu}}\cdot \wi\mu&\hbox{with}& Z_{\wi\mu}\,\df\, \wi\mu[\psi]\\
\nu\,\df\, \frac{\psi}{Z_{\wi\nu}}\cdot \wi\nu&\hbox{with}& Z_{\wi\nu}\,\df\, \wi\nu[\psi]\eq
Let $\wi f$ and $f$ stand for the Radon-Nikodym densities of $\wi \mu$ with respect to $\wi\nu$
and of $\mu$ with respect to $\nu$. Obviously, we have
\bq
f&=&\frac{Z_{\wi\nu}}{Z_{\wi\mu}}\wi f\eq
Finally, choose $\wi m$ and $m$ to be  medians of $\wi f$ and $f$ with respect to $\wi\nu$ and $\nu$.
The following result is well-known.
\begin{lem}
We have\bq
\int \lve f-m\rve\,d\nu\,\leq \, \lVe \mu-\nu\rVe_{\mathrm{tv}}\,\leq\, 2\int \lve f-m\rve\,d\nu\\
\int \lve \wi f-\wi m\rve\,d\wi\nu\,\leq \, \lVe \wi\mu-\wi\nu\rVe_{\mathrm{tv}}\,\leq\, 2\int \lve \wi f-\wi m\rve\,d\wi\nu
\eq
\end{lem}
\proof
Of course it is sufficient to show the bounds for $ \lVe \mu-\nu\rVe_{\mathrm{tv}}$. 
They are a consequence of 
\bq
\lVe \mu-\nu\rVe_{\mathrm{tv}}&=&\int \lve f-1\rve\,d\nu\eq
and of the following characterization of a median:
\bqn{median}
\int \lve f-m\rve\,d\nu&=&\inf\lt\{\int \lve f-r\rve\,d\nu\st r\in\RR\rt\}\eqn
So the lower bound is immediate and for the upper bound, just note that
\bq
\lve 1-m\rve&=&\lve \int (f-m)\,d\nu\rve\\
&\leq &  \int \lve f-m\rve\,d\nu\eq
\wwtbp
The interest of the introduction of the medians comes from:
\begin{lem}\label{essinf}
We have
\bq
\int \lve f-m\rve\,d\nu&\leq &  \frac{\psi_\wedge}{\psi_\vee}\int \lve \wi f-\wi m\rve\,d\wi\nu
\eq
and it follows from the previous lemma that 
\bq
  \frac{\psi_\vee}{2\psi_\wedge}\lVe \wi\mu-\wi\nu\rVe_{\mathrm{tv}}\ \leq \ 
 \lVe \mu-\nu\rVe_{\mathrm{tv}}\ \leq \ 2 \frac{\psi_\wedge}{\psi_\vee}\lVe \wi\mu-\wi\nu\rVe_{\mathrm{tv}}\eq
\end{lem}
\proof
From \eqref{median}, we have
\bq
\int \lve \wi f-\wi m\rve\,d\wi\nu&=&
\inf\lt\{\int \lve \wi f-r\rve\,d\wi \nu\st r\in\RR\rt\}\\
&=&\inf\lt\{\int \lve \frac{Z_{\wi\nu}}{Z_{\wi\mu}} f-r\rve\,d\wi \nu\st r\in\RR\rt\}\\
&=&\frac{Z_{\wi\nu}}{Z_{\wi\mu}} \inf\lt\{\int \lve f-r\rve\frac{\psi}{Z_{\wi\nu}} \,d \nu\st r\in\RR\rt\}\\
&\geq &\frac{Z_{\wi\nu}}{Z_{\wi\mu}}\frac{\essinf_{\wi\nu}\psi}{Z_{\wi\nu}} \inf\lt\{\int \lve f-r\rve\,d\nu\st r\in\RR\rt\}\\
&=& \frac{\essinf_{\wi\nu}\psi}{Z_{\wi\mu}}\int \lve f-m\rve\,d\nu\\
&\geq & \frac{\essinf_{\wi\nu}\psi}{\esssup_{\wi\mu}\psi}\int \lve f-m\rve\,d\nu\\
&\geq &
 \frac{\psi_\wedge}{\psi_\vee}\int \lve f-m\rve\,d\nu
\eq
\wwtbp\par
For any fixed $t\geq 0$, it remains to apply these general bounds with
\bq 
\wi \mu&\df& \wi\mu_0 \wi P_t\\
\wi\nu&\df& \wi\eta\\
\psi&\df& 1/\varphi
\eq
\par
Since  $ \psi_\wedge/\psi_\vee= \varphi_\vee/\varphi_\wedge$,
the conclusion of Lemma~\ref{essinf} implies the wanted bound.
\begin{rem}
From \eqref{mutwimu0}, we could have been tempted to write that for any $f\in\cF$,
\bq
\mu_t[f]-\nu[f]&=&
\frac{\wi\mu_0[\wi P_t[f/\varphi]]}{\wi\mu_0[\wi P_t[1/\varphi]]}-\frac{\wi\eta[f/\varphi]}{\wi \eta[1/\varphi]}\\
&=&
\frac1{\wi\mu_0[\wi P_t[1/\varphi]]}(\wi\mu_0[\wi P_t[f/\varphi]]-\wi\eta[f/\varphi])+\frac{\wi\eta[f/\varphi]}{\wi \eta[1/\varphi]\wi\mu_0[\wi P_t[1/\varphi]] }(\wi \eta[1/\varphi]-\wi\mu_0[\wi P_t[1/\varphi]])\\
&\leq & \varphi_\vee\lve \wi\mu_0[\wi P_t[f/\varphi]]-\wi\eta[f/\varphi]\rve+ \frac{ \varphi_\vee^2}{\varphi_\wedge}\lVe f\rVe_\iy \lve \wi \eta[1/\varphi]-\wi\mu_0[\wi P_t[1/\varphi]])\rve
\eq
Taking the supremum of $f\in\cF$ satisfying $\lVe f\rVe_\iy\leq 1$, it appears that
\bq
\lVe \mu_t-\nu\rVe_{\mathrm{tv}}&\leq & \lt(\lt(\frac{\varphi_\vee}{\varphi_\wedge}\rt)+\lt(\frac{\varphi_\vee}{\varphi_\wedge}\rt)^2\rt)\lVe \wi\mu_0\wi P_t -\wi\eta\rVe_{\mathrm{tv}}
\\&\leq & 2\lt(\frac{\varphi_\vee}{\varphi_\wedge}\rt)^2\lVe \wi\mu_0\wi P_t -\wi\eta\rVe_{\mathrm{tv}}\eq
which is worse than the bound of Theorem \ref{tvtv} by a factor ${\varphi_\vee}/{\varphi_\wedge}$.
\end{rem}

\subsection{Proof of Theorem \ref{tvL2}}

Since Theorem \ref{tvtv} brings us back to the situation of convergence to equilibrium of Markov processes, it is sufficient to use the argument 
of Fill \cite{MR1097464}) for non-reversible processes. We recall them below for the sake of completeness.
\par\sm
To gain a factor 2, it is in fact better not to use Theorem \ref{tvtv}, but to directly make a comparison between $\LL^2$ quantities.
More precisely,
for given $\mu_0\in\cP$ and $t\geq 0$, denote by $f_t$ (respectively $\wi f_t$) the density of the probability 
$\mu_t$  with respect to $\nu$  (resp.\ 
$\wi\mu_t\df \wi\mu_0\wi P_t$ with respect to $\wi \eta$).
We have by the Cauchy-Schwarz inequality,
\bq
\lVe \mu_t-\nu\rVe_{\mathrm{tv}}&=&\sum_{x\in S}\lve f_t(x)-1\rve\nu(x)\\
&\leq & \sqrt{I_t}\eq
where
\bq
I_t&\df& \sum_{x\in S}( f_t(x)-1)^2\,\nu(x)\eq
Let us also define 
\bq
\wi I_t&\df& \sum_{x\in S}(\wi f_t(x)-1)^2\,\wi\eta(x)\eq
It is easy to compare these quantities:
\begin{lem}\label{II}
For any $t\geq 0$, we have
\bq
I_t&\leq &\lt( \frac{\varphi_\vee}{\varphi_\wedge}\rt)^2\wi I_t\eq
\end{lem}
\proof
One recognizes in $I_t$ the variance of $f_t$ with respect to $\nu$, so that
\bq
I_t&=&\inf\lt\{ \sum_{x\in S} (f_t(x)-r)^2\, \nu(x)\st r\in\RR\rt\}\eq
Similarly we have
\bq\wi I_t&=&\inf\lt\{ \sum_{x\in S} (\wi f_t(x)-r)^2\, \wi\eta(x)\st r\in\RR\rt\}\eq
The same arguments as those used 
 in Lemma \ref{essinf} give the conclusion without difficulty.\wwtbp
 Putting together these estimates, we end up with
 \bq
\lVe \mu_t-\nu\rVe_{\mathrm{tv}}&\leq & 
\frac{\varphi_\vee}{\varphi_\wedge}\sqrt{\wi I_t}\eq
To study the evolution of $\wi I_t$ with respect to the time $t\geq 0$, recall that
\bq
\fo x\in S,\,\fo t\geq 0,\qquad \pa_t \wi f_t(x)&=&\wi L^*[\wi f_t](x)\eq
This comes from the relation $\wi f_t=\wi P_t^*[f_0]$, where $\wi P_t^*$ the adjoint operator of $\wi P_t$ in $\LL^2(\wi\eta)$.
Thus we get that for all $t\geq 0$,
\bqn{LLdiamond}
\nonumber\pa_t \wi I_t&=&2 \wi\eta[ ( \wi f_t-1)\pa_t \wi f_t]\\
\nonumber&=&2\wi\eta[  ( \wi f_t-1)\wi L^*[\wi f_t]]\\
\nonumber&=&2\wi\eta[  ( \wi f_t-1)\wi L^*[\wi f_t-1]]\\
\nonumber&=&2\wi\eta[ \wi L[ \wi f_t-1](\wi f_t-1)]\\
&=&2\wi\eta[\wi  L^\diamond [ \wi f_t-1](\wi f_t-1)]\eqn
By definition of $\wi\lambda$,  the r.h.s.\ is bounded
above by $-2\wi\lambda \wi I_t$, which leads to the ordinary differential inequality
\bq
\fo t\geq 0,\qquad\pa_t \wi I_t&\leq &-2\wi\lambda \wi I_t\eq
Gronwall's lemma implies that
\bq
\fo t\geq 0,\qquad \wi I_t&\leq & \exp(-2\wi\lambda t)\wi I_0\eq
so it remains to bound $\wi I_0$ above.
But note that
\bq
\wi I_0&=& \wi\eta[\wi f_0^2]-1\\
&\leq &  \wi\eta[\wi f_0^2]\\
&=&\wi\mu_0[\wi f_0]\\
&\leq &\lVe\wi f_0 \rVe_\iy\\
&\leq &\frac1{\wi\eta_\wedge}\\
&=&
\frac{\eta[\varphi\varphi^*]}{(\varphi\varphi^*\eta)_{\wedge}}
\eq
which, in conjunction with Theorem \ref{tvtv}, leads to the bound of Theorem \ref{tvL2}.\par
The expression \eqref{Ldiamond} for $\wi L^\diamond$ is a consequence of
\bqn{wiLst}
\nonumber\fo x\not=y,\qquad \wi L^*(x,y)&=&\frac{\wi\eta (y)}{\wi\eta(x)}\wi L(y,x)\\
&=&\frac{\varphi^*(y)\eta(y)}{\varphi^*(x)\eta(x)}L(y,x)\eqn
\par
The Poincaré formulation \eqref{Poincare} comes from the variational characterization of the eigenvalues and from the equality 
\bq \fo g\in\cF,\qquad \wi\eta[g\wi L^\diamond[g]]&=&\wi\eta[g\wi L[g]]\eq
already used in \eqref{LLdiamond}.
\par
\begin{rem}\label{rates2}
Similarly to the lower bound in Theorem \ref{tvtv}, we have also in Lemma \ref{II}
\bq
\fo t\geq 0,\qquad
I_t&\geq &\lt( \frac{\varphi_\wedge}{\varphi_\vee}\rt)^2\wi I_t\eq
In particular $\sqrt{I_t}$ and $\sqrt{\wi I_t}$  have the same asymptotic exponential rate of convergence. This common rate is
the smallest real part of the non-zero eigenvalues of $-\wi L$, but since this operator is not assumed to be reversible,
this rate may be larger than $\wi\lambda$.
\end{rem}
\par
\begin{rem}
It is possible to improve the pre-exponential factor $\sqrt{\frac{\eta[\varphi\varphi^*]}{(\varphi\varphi^*\eta)_{\wedge}}}
\frac{\varphi_\vee}{\varphi_\wedge}$ in Theorem \ref{tvL2}, but  at the expense of the rate $\wi\lambda$, via the logarithmic Sobolev inequalities associated to
the symmetrization $\wi L^\diamond$ of $\wi L$.
\par
Let $\wi\alpha>0$ be the largest constant such that for all $g\in\cF$,
\bqn{lSi}
\wi\alpha \sum_{x\in S}g^2(x)\ln\lt(\frac{g^2(x)}{\wi\eta[g^2]}\rt)\, \varphi^*(x)\varphi(x)\eta(x)&\leq & \sum_{x,y\in S}(g(y)-g(x))^2\, \varphi^*(x)\varphi(y)\eta(x) L(x,y)\eqn
\par
Then we have
\bqn{ent}
\sup_{\mu_0\in\cP}\lVe \mu_t-\nu\rVe_{\mathrm{tv}}&\leq &
\sqrt{2\ln\lt(\frac{\eta[\varphi\varphi^*]}{(\varphi\varphi^*\eta)_{\wedge}}\rt)
\frac{\varphi_\vee}{\varphi_\wedge}}\exp(-(\wi\alpha/2) t)\eqn
\par
The proof of this bound has the same structure as the one of Theorem \ref{tvL2}, with the quantities $I_t$ and $\wi I_t$
replaced by the relative entropies
\bq
J_t&\df& \sum_{x\in S} f_t(x)\ln(f_t(x))\, \nu(x)\\
\wi J_t&\df& \sum_{x\in S} \wi f_t(x)\ln(\wi f_t(x))\, \wi\eta(x)\eq
\par
Indeed, Pinsker's inequality gives the bound
\bq
\lVe \mu_t-\nu\rVe_{\mathrm{tv}}&\leq & \sqrt{2}\sqrt{J_t}\eq
Next, taking into account the relations (see Holley and Stroock \cite{MR90g:60091})
\bq
J_t&=&\inf\lt\{ \sum_{x\in S}(f_t(x)\ln(f_t(x)) -f_t(x)\ln(r)-f_t(x)+r)\,\nu(x)\st r\in\RR_+\rt\}\\
\wi J_t&=&\inf\lt\{ \sum_{x\in S}(\wi f_t(x)\ln(\wi f_t(x)) -\wi f_t(x)\ln(r)-\wi f_t(x)+r)\,\wi\eta(x)\st r\in\RR_+\rt\}
\eq
we deduce as in Lemma \ref{essinf} and Remark \ref{rates1} that
\bq
\fo t\geq 0,\qquad \frac{\varphi_\wedge}{\varphi_\vee}\wi J_t \ \leq \ J_t\ \leq \ \frac{\varphi_\vee}{\varphi_\wedge}\wi J_t\eq\par
As a consequence, we get
 \bq
\fo t\geq 0,\qquad \lVe \mu_t-\nu\rVe_{\mathrm{tv}}&\leq & 
\sqrt{2\frac{\varphi_\vee}{\varphi_\wedge}}\sqrt{\wi J_t}\eq
which reduces our task to the investigation of the time evolution of $\wi J_t$.
\par
By differentiation, it appears that
\bq
\pa_t \wi J_t&=&\sum_{x\in S} (1+\ln(\wi f_t(x))) \pa_t \wi f_t(x)\, \wi\eta(x)\\
&=&\sum_{x\in S} (1+\ln(\wi f_t(x))) \wi L^*[\wi f_t(x)]\, \wi\eta(x)\\
&=&\sum_{x,y\in S} (1+\ln(\wi f_t(x))) (\wi f_t(y)-\wi f_t(x))\, \wi\eta(x)\wi L^*(x,y)\eq
To proceed, note (cf.\ for instance Miclo \cite{MR1478724}) that for all $ x,y\in S$,
\bq
(1+\ln(\wi f_t(x))) (\wi f_t(y)-\wi f_t(x))&\leq & \wi f_t(y)\ln( \wi f_t(y))- \wi f_t(x)\ln( \wi f_t(x))-\lt(\sqrt{\wi f_t(y)}-\sqrt{\wi f_t(x)}\rt)^2\eq
and that by invariance of $\wi\eta$ with respect to $\wi L^*$,
\bq\sum_{x,y\in S} \lt(\wi f_t(y)\ln( \wi f_t(y))- \wi f_t(x)\ln( \wi f_t(x))\rt)\, \wi\eta(x)\wi L^*(x,y)&=&0\eq
Thus we end up with the refined Jensen type bound:
\bq\pa_t \wi J_t&\leq &-\sum_{x,y\in S} \lt(\sqrt{\wi f_t(y)}-\sqrt{\wi f_t(x)}\rt)^2\, \wi\eta(x)\wi L^*(x,y)\\
&=&-\sum_{x,y\in S} \lt(\sqrt{\wi f_t(y)}-\sqrt{\wi f_t(x)}\rt)^2\, \wi\eta(y)\wi L^*(y,x)\\
&=&-\sum_{x,y\in S} \lt(\sqrt{\wi f_t(y)}-\sqrt{\wi f_t(x)}\rt)^2\, \varphi^*(x)\varphi(y)\eta(x) L(x,y)\eq
where we used \eqref{wiLst}.
The logarithmic Sobolev inequality \eqref{lSi}, with $g\df f_t$, allows  comparison of the r.h.s.\ with $\wi J_t$ to give the
differential inequality
\bq
\fo t\geq 0,\qquad\pa_t \wi J_t&\leq &-\wi\alpha \wi J_t\eq
Gronwall's lemma implies again that
\bq
\fo t\geq 0,\qquad \wi J_t&\leq & \exp(-\wi\alpha t)\wi J_0\\
&\leq & \exp(-\wi\alpha t)\ln( (\wi f_0)_\vee)\\
&\leq & \exp(-\wi\alpha t)\ln(1/\wi\eta_\wedge)
\eq
The announced bound \eqref{ent} follows.\par\sm
Despite the deterioration of exponential rate in \eqref{ent},
this bound can be interesting for not too large times $t\geq 0$, especially when one looks for ``quasi-mixing times".
Diaconis and Saloff-Coste \cite{MR1410112} have shown the following general bound between the logarithmic Sobolev constant $\wi\alpha$ and the spectral gap $\wi\lambda$:
\bqn{alla}
\wi\alpha&\geq & \frac{1-2\wi\eta_\wedge}{\ln(1/\wi\eta_\wedge-1)}\wi\lambda\eqn
(where the factor on the r.h.s.\ is taken to be $1/2$ in the particular case where $\wi\eta_\wedge=1/2$). 
But this relation is not very pertinent for quasi-mixing times estimates: if $\tau_{\wi\lambda}\geq 0$ and $\tau_{\wi\alpha}\geq 0$
are  the times $t\geq 0$ in Theorem \ref{tvL2} and \eqref{ent} such that the corresponding upper bounds are equal to 1,
we get
\bq
\tau_{\wi\lambda}&=&\frac{1}{\wi\lambda}(\ln(\varphi_{\vee}/\varphi_\wedge)+\ln(1/ \wi\eta_\wedge))\\
\tau_{\wi\alpha}&=&\frac{2}{\wi\alpha}(\ln(\varphi_{\vee}/\varphi_\wedge)+\ln(\ln(1/ \wi\eta_\wedge)))\eq
and the injection of \eqref{alla}  leads to the disappointing $\tau_{\wi\lambda}\ll\tau_{\wi\alpha}$ for small $\wi\eta_\wedge>0$.
Indeed, the interest of \eqref{ent} appears when one has good estimates on $\wi\alpha$ (by tensorization for instance)
and $\wi\eta_\wedge$ is very small. Simple examples on product spaces are provided in Subsection \ref{ape}.
Nevertheless, we believe that modified logarithmic Sobolev inequalities (see e.g.\ the article of Bobkov and Tetali \cite{MR2283379}), namely the consideration
of the best constant $\wit\alpha>0$ such that for all $g\in\cF$,
\bq
\lefteqn{\wit\alpha \sum_{x\in S}g^2(x)\ln\lt(\frac{g^2(x)}{\wi\eta[g^2]}\rt)\, \varphi^*(x)\varphi(x)\eta(x)}\\
&\leq & \sum_{x,y\in S}(\vert g(y)\vert -\vert g(x)\vert )(\ln(\vert g(y)\vert )-\ln(\vert g(x)\vert ))\, \varphi^*(x)\varphi(y)\eta(x) L(x,y)\eq
is  better suited to the above entropic approach.
\end{rem}

\subsection{Proof of Theorem \ref{tvL2rev}}

Under the assumption that $\nu$ is reversible for $L$, we have that $L^*=L$. The equations for $\varphi$ and $\varphi^*$
are thus the same and only the corresponding renormalizations are different. If follows that $\varphi$ and $\varphi^*$
are proportional and since only ratios enter the pre-exponential factor of Theorem \ref{tvL2},
it can be replaced by the pre-exponential factor of Theorem \ref{tvL2rev} (recall the normalization $\eta[\varphi^2]=1$).
\par
But the main advantage of Theorem \ref{tvL2rev} is the explicit rate $\lambda_2-\lambda_1$.
It is a consequence of the conjugacy relation \eqref{conj}. It shows first that $\wi L$ must  be reversible with respect 
to $\wi\eta$ (but this can also be checked directly from the expressions \eqref{wiL} and \eqref{wieta}) and second
that the spectrum of $ \wi L$ is obtained from the spectrum of $L-V$ by subtracting the value $\lambda_1$.
In particular the spectral gap $\wi\lambda$ of $\wi L^\diamond=\wi L$ is equal to $\lambda_2-\lambda_1$.
\begin{rems}
(a) The fact that the spectrum of $ \wi L$ is obtained from the spectrum of $L-V$ by subtracting the value  $\lambda_1$
is always true, but in the non-reversible case it is not clear how to use this possibly complex valued spectrum
to deduce a bound on $\wi\lambda$.
In the reversible situation Remark \ref{rates2} can be made more precise: the common asymptotic exponential rate of 
$\sqrt{I_t}$ and $\sqrt{\wi I_t}$ is $\lambda_2-\lambda_1$.
\par\sm
(b) The logarithmic Sobolev inequality approach is equally valid in the reversible case, we get
\bq
\sup_{\mu_0\in\cP}\lVe \mu_t-\nu\rVe_{\mathrm{tv}}&\leq &
\sqrt{2\ln\lt(\frac{1}{(\varphi^2\eta)_{\wedge}}\rt)
\frac{\varphi_\vee}{\varphi_\wedge}}\exp(-(\wi\alpha/2) t)\eq
where $\wi\alpha$ is the logarithmic Sobolev constant associated to the symmetric operator $\wi L$ in $\wi\eta$
(in particular \eqref{alla} is satisfied with $\wi\lambda$ replaced by $\lambda_2-\lambda_1)$.\nopagebreak
\end{rems}

\section{Examples}\label{examples}

Several basic examples are provided here, which in particular serve to illustrate
some assertions made in the previous theoretical developments.

\subsection{A birth and death example with $\lambda_1\approx \lambda_2-\lambda_1$}\label{abadewlsll}

This example and the next two are birth and death processes on $\bar S\df \lin 0,N\rin$, with $N\in\NN$,
absorbed in $0$. So $S=\lin 1,N\rin$, $\iy=0$ and  $L$ gives positive rates only to the oriented edges 
$(x,x+1)$ and $(x+1,x)$ where $x\in \lin 1,N-1\rin$.
In this one-dimensional setting, $L$ admits a unique reversible probability $\eta$.
Let us assume that the killing rate in $1$ is 1, namely $V(1)=\bar L(1,0)=1$.
The other values of $V$ are taken to be zero.\par\sm
Specifically for this example, we choose
\bqn{L1}
\fo x\in \lin 1,N-2\rin,& & L(x,x+1)\ \df\ L(x+1,x)\ \df\ 1\\
\label{L1b}&&\hskip-8mm L(N-1,N)\ =\ 1\ \hbox{ and }\ L(N,N-1)\ =\ 2
\eqn
(the value 2 simplifies the analysis of the reflection at $N$ by replacing the forbidden jump to $N+1$ by a supplementary
jump at $N-1$).
The reversible probability $\eta$ is then given by
\bq
\fo x\in S,\qquad \eta(x)&=&\frac1N\eq
Let $\varphi$ be the function defined by
\bqn{ex1} 
\fo x\in S,\qquad \varphi(x)&\df& \frac1Z\sin(\pi x/(2N))\eqn
where $Z$ is the renormalization constant such that $\eta[\varphi^2]=1$.
Due to the value 2 in \eqref{L1b}, it is easy to check that $(L-V)[\varphi]=2(\cos(\pi/(2N))-1)\varphi$.
The positivity of $\varphi$ and Perron-Frobenius theorem imply that $\varphi$ is
indeed the function considered in the introduction and that
\bq
\lambda_1&=&2(1-\cos(\pi/(2N)))\eq
The density of the quasi-invariant probability measure $\nu$ with respect to $\eta$ is proportional to $\varphi$.
\par
More generally, define for $k\in \lin 1,N-1\rin$, the function $\varphi_k$ by
\bq
\fo x\in S,\qquad \varphi_k(x)&\df&\sin((2k+1)\pi x/(2N))\eq
By straightforward calculation, $(L-V)[\varphi_k]=2(\cos((2k+1)\pi/(2N))-1)\varphi_k$. Thus the spectrum of $L-V$ is
$\{2(\cos((2k+1)\pi/(2N))-1)\st k\in \lin 0,N-1\rin\}$. In particular
\bq
\lambda_2-\lambda_1&=&2(\cos(\pi/(2N))-\cos(3\pi/(2N))\\
&=&2\sin(\pi/N)\sin(\pi/(2N))\\
&=& \frac{\pi^2}{N^2}(1+\cO(N^{-2}))\eq
as $N$ goes to infinity. Since \bq\lambda_1&=&\frac{\pi^2}{4N^2}(1+\cO(N^{-2}))\eq
in this situation $\lambda_1$ and $\lambda_2-\lambda_1$ are of the same order, meaning that absorption and convergence to quasi-stationarity
happen at similar rates.\par
From \eqref{ex1}, we deduce that
\bq
\frac{\varphi_\vee}{\varphi_\wedge}&=&\frac1{\sin(\pi /(2N))}\\
&= & \frac{2N}{\pi}(1+\cO(N^{-2}))\eq
Taking into account the classical Riemann sum approximation, we furthermore get
\bq
Z^2&=&\frac1N\sum_{x\in\lin 1, N\rin}\sin^2(\pi x/(2N))\\
&=&(1+\cO(N^{-1})) \int_0^1\sin^2(\pi u/2)\,du\\
&=&\frac12(1+\cO(N^{-1}))\eq
The first bound of Theorem \ref{tvL2rev} asserts that
\bq
\sup_{\mu_0\in\cP}\lVe \mu_t-\nu\rVe_{\mathrm{tv}}&\leq &
\frac{2\sqrt{2}}{\pi^2}N^{5/2}\exp\lt(- \frac{\pi^2}{N^2}t(1+\cO(N^{-2}))\rt)(1+\cO(N^{-1}))\eq
(the second bound of Theorem \ref{tvL2rev}, which doesn't need the estimate on $Z$, leads to a similar bound with $2\sqrt{2}$ replaced by 4).
It follows that for any given $s>0$, if 
\bq
t&=&\frac5{2\pi^2}N^2\ln(N)+\frac{s}{\pi^2}N^2\eq
then 
\bq
\sup_{\mu_0\in\cP}\lVe \mu_t-\nu\rVe_{\mathrm{tv}}&\leq &\frac{2\sqrt{2}}{\pi^2}(1+\cO(N^{-1}))\exp(-s)\eq

\subsection{A birth and death example with $\lambda_1\ll \lambda_2-\lambda_1$}\label{abadewlll}

The setting is as in the previous example, except that for some $r>1$, we replace \eqref{L1} and \eqref{L1b} by
\bqn{L2}
\fo x\in \lin 1,N-1\rin,&\ & \lt\{
\begin{array}{ccc}L(x,x+1)& \df& r\\ L(x+1,x)& \df&1
\end{array}\rt.\\
\label{L2b}&&\hskip-20mm L(N-1,N)\ =\ 1\ \hbox{ and }\ L(N,N-1)\ =\ 1+r
\eqn
The reversible probability $\eta$ is then given by
\bqn{etar}
\fo x\in S,\qquad \eta(x)&=&\frac{r-1}{r^{N}-1}r^{x-1}\eqn
Contrary to the previous example, it seems more difficult to derive explicit formulas for the eigenvalues and eigenfunctions associated to $L-V$.
To describe them, consider the rational fraction in $X$,\bq
P_N(X)&\df& \frac{X^{2(N+1)}-X^{2N}+r^{1-N}X^2-r^{-N-1}}{X^2-r^{-1}}\eq
\begin{lem}\label{Lambda}
For $N\geq 1$, $P_N$ is a polynomial which admits $2N$ distinct zeros. Denote by $R$ the set of zeros.
Let $\Lambda$ be the image of $R$ by the mapping 
\bq
\Psi\st\rho \mapsto \frac{(1+r)\rho-1-r\rho^2}{\rho}\eq
For $N> (1+r)/(r-1)$, 
the spectrum of $V-L$ is $\Lambda$ and for any $\lambda\in\Lambda$,
an associated eigenfunction $\varphi_\lambda$ is defined by
\bq
\fo x\in S,\qquad \varphi_\lambda(x)&\df& \rho_+^x-\rho_-^x\eq
where
\bqn{rhopm}
\rho_{\pm}&\df&\frac{1}{2r}(r+1-\lambda\pm\sqrt{(\lambda-1-r)^2-4r})\eqn
(with $\sqrt{\cdot}$ standing for the principal value of the complex square root) are the reciprocal images of $\lambda$ by $\Psi$.\end{lem}
\proof
Let $\lambda$ be an eigenvalue of $V-L$ and $\varphi$ be an associated eigenfunction on $S$.
With the convention that $\varphi(0)=0$, the values of $\varphi$ satisfy the
recursive formula
\bqn{recur}
\fo x\in \lin 1,N-1\rin,\qquad
\varphi(x+1)&=&\frac{(1+r-\lambda)\varphi(x)-\varphi(x-1)}{r}\eqn
It follows that on $\lin 1, N\rin$, $\varphi$ is necessarily proportional to the functions $\varphi_\lambda$
defined above, where $\rho_\pm$ are the solutions of the quadratic equation in $X$,
\bqn{quad}
rX^2+(\lambda-1-r)X+1&=&0\eqn
except if this equation admits a double solution $\rho_*$, in which case $\varphi$ must be proportional to the function $\varphi_*$ defined by 
\bq
\fo x\in S,\qquad \varphi_*(x)&\df& x\rho_*^x\eq
Whatever the case,  we have that
\bq
\fo x\in\lin 1,N-1\rin,\qquad (L-V)[\varphi](x)&=&-\lambda \varphi(x)\eq
This relation is also satisfied at $x=N$ if and only if $\varphi(N+1)=\varphi(N-1)$
(as in the previous example, this justifies the simplifying choice of $L(N,N-1)=1+r$).\par\sm
$\bullet$ Let us first consider the situation where \eqref{quad} admits a double solution.
One computes immediately that this corresponds to $\lambda=(1\pm \sqrt{r})^2$ and $\rho_*=\mp 1/\sqrt{r}$.
The condition $\varphi_*(N+1)=\varphi_*(N-1)$ is equivalent to 
$\rho_*=\pm\sqrt{(N-1)/(N+1)}$. The assumption $N> (1+r)/(r-1)$ forbids that
$\sqrt{(N+1)/(N-1)}=\sqrt{r}$, so that we are led to a contradiction. Only the next case is possible.\par\sm
$\bullet$ 
Assume that \eqref{quad} has two distinct solutions $\rho_+$ and $\rho_-$, they are  given in the statement of the above lemma.
The condition $\varphi(N+1)=\varphi(N-1)$ amounts to 
\bq
\rho_+^{N+1}-\rho_+^{N-1}-\rho_-^{N+1}+\rho_-^{N-1}&=&0\eq
But from \eqref{quad} we see that $\rho_-=1/(r\rho_+)$, so $\rho_+$ is a solution of
\bq
X^{2(N+1)}-X^{2N}+r^{1-N}X^2-r^{-N-1}&=&0\eq
This equation admits two obvious solutions, $X=1/\sqrt{r}$ and $X=-1/\sqrt{r}$,
so that $P_N$ is indeed a polynomial.
But these values are not allowable for $\rho_+$, because  we would have $\rho_+=\rho_-$.
It follows that $\rho_+$ is a root of $P_N$.
Note that if $\rho\in\CC$ is a root of $P_N$, the same is true for $1/(r\rho)$ and that
$1/\sqrt{r}$ and $-1/\sqrt{r}$ are the only fixed points of the involutive mapping $\xi\st \CC\setminus\{0\}\ni\rho\mapsto 1/(r\rho)$.
As a consequence, we can group the roots of $P_N$ by pairs stable by $\xi$, say
$\{\rho_1,\xi(\rho_1)\}$, $\{\rho_2,\xi(\rho_2)\}$, ..., $\{\rho_N,\xi(\rho_N)\}$.
Moreover, notice that the mapping $\Psi$ defined in the above lemma is constant on each of these pairs,
so that the cardinality of $\Lambda\df \Psi(R)$ is at most $N$.
But it appears from \eqref{quad} and from the previous discussion that
all the eigenvalues of $V-L$ are elements of $\Lambda$.
From \eqref{recur} we deduce that all the eigenvalues of $V-L$ are simple
and since by reversibility $V-L$ is known to be diagonalizable, it follows
$V-L$ admits $N$ distinct eigenvalues. Thus $\Lambda$ must be of cardinality $N$
and exactly consists of the eigenvalues of $V-L$.
It is interesting to remark that it is relatively difficult to check directly 
that all the $\rho_1$, $\rho_2$, ..., $\rho_N$ are different, or equivalently that
 all the roots of $P_N$ are distinct (try to compute its discriminant).\wwtbp
\par
By the Perron-Frobenius theorem, the smallest eigenvalue $\lambda_1$ of $V-L$
is characterized by the fact that the associated eigenfunction
has a fixed sign. This observation in conjunction with
Lemma \ref{Lambda} lead to
\begin{pro} \label{rgeq1}
For large $N$ we have
\bq
\lambda_1\sim \frac12(r+1)(r-1)^2\frac1{r^{N+1}}\eq
moreover, if $\varphi$ is an associated eigenvector, 
\bq
\frac{\varphi_\vee}{\varphi_\wedge}&=& \frac{r}{r-1}(1+\cO(r^{-N}))\eq
\end{pro}
\proof
With the notation of Lemma \ref{Lambda}, we have $\rho_+\rho_-=1/r$ (recall \eqref{quad}).
So if $\rho_+>0$, then we get $0<\rho_-<\rho_+$. It follows that the mapping
$\RR_+^*\ni u\mapsto \rho_+^u -\rho_-^u$ does not vanish and in particular
$\varphi_\lambda$ only takes positive values. By consequence, the corresponding $\lambda\in\Lambda$
is the first Dirichlet eigenvalue $\lambda_1$.
To work out this program, we begin by showing that for $N$ large enough, there exists $\rho_1\in R\cap\RR_+$ satisfying 
\bqn{rho1}
\frac1{r^2}\ \leq \ \rho_1^2\ \leq \ \frac1{r^2}+\frac1{r^{N+1}}\eqn
 It is enough to show that there exists
 $\rho_1\in [1/r,\sqrt{1/r^2+r^{-1-N}}]$ such that
\bqn{rho1b} Q(\rho_1^2)\ \df\  \rho_1^{2(N+1)}-\rho_1^{2N}+r^{1-N}\rho_1^2-r^{-N-1}\ =\ 0\eqn
 Write $h_1=r^{N+1}(\rho_1^2-1/r^2)$ and for all $h\geq 0$,
 \bq
 f(h)&\df& r^{2N}Q\lt(\frac1{r^2}+\frac{h}{r^{N+1}}\rt)\\&=&
 \lt(1+\frac{h}{r^{N-1}}\rt)^N\lt(\frac1{r^2}-1+\frac{h}{r^{N+1}}\rt)+h\eq
 we just need to check that $f(0)\leq 0$ and $f(1)\geq 0$.
 The former inequality is immediate and the latter one is satisfied for $N$ large enough, since
 $\lim_{N\ri\iy} f(1)=1/r^2$.\\
Next, injecting the a priori bound \eqref{rho1} in \eqref{rho1b}, it follows that
\bq
\lim_{N\ri\iy}h_1
&=&1-\frac1{r^2}\eq
Replacing $\rho_1=\sqrt{r^{-2}+(1-r^{-2})r^{-(N+1)}(1+\circ(1))}=r^{-1}+(1-r^{-2})r^{-N}(1/2+\circ(1))$ in
\bq\lambda_1&=&
 \frac{(1+r)\rho_1-1-r\rho_1^2}{\rho_1}\eq
 we deduce the first announced behavior.
Let $\rho_{1-}$ and $\rho_{1+}$ be the corresponding values of $\rho_-$ and $\rho_+$, 
from $\rho_{1-}\rho_{1+}=1/r$, we obtain 
\bqn{rho1-}
\rho_{1-}\ =\ \rho_1\ =\ \frac1{r}+\cO(r^{-N})
&\hbox{ and }&   \rho_{1+}\ =\ 1+\cO(r^{-N}) \eqn
Taking into account the expression of $\varphi\df\varphi_{\lambda_1}$ given in Lemma \ref{Lambda},
we get
\bq
\frac{\varphi_\vee}{\varphi_\wedge}&=&\frac{\varphi(N)}{\varphi(1)}\\
 &=&\frac{ \rho_{1+}^N-\rho_{1-}^N}{ \rho_{1+}-\rho_{1-}}\\
 &=&\frac{r}{r-1}(1+\cO(r^{-N}))\eq
 \wwtbp
 \par
 To be in position to use Theorem \ref{tvL2rev}, it remains to evaluate $\lambda_2-\lambda_1$.
 \par
 From the previous proof, it appears there is only one eigenvalue $\lambda\in\Lambda$
 such that $\rho_->0$. Moreover there is at most one eigenvalue $\lambda\in\Lambda$
 such that $\rho_-<0$. Indeed, in this case we have $\rho_-<\rho_+<0$
 and it follows  from Lemma \ref{Lambda} that
 $\varphi_{\lambda}(x)>0$ for $x\in S$ odd and $\varphi_{\lambda}(x)<0$ for $x\in S$ even,
 in particular $\varphi_{\lambda}$ has the maximal number of sign changes.
 The discrete version of Sturm's theorem (see for instance Miclo \cite{MR2438701}) then
 implies that $\lambda$ must be $\lambda_N$, the largest eigenvalue of $V-L$.
 Since $R$ is symmetrical with respect to zero, $-\rho_{1-}$ and $-\rho_{1+}$
 (with the notation of \eqref{rho1-}) also belong to $R$ and this 
 leads to the estimate 
 \bq
 \lambda_N&=& 2(1+r)+\cO(r^{-N})\eq
 \par
 The previous arguments show that except for $\rho_{1-}$, $\rho_{1+}$, $-\rho_{1-}$ and $-\rho_{1+}$,
all the other elements of $R$ are complex numbers which are not real.
It follows from Lemma  \ref{Lambda} that for $\lambda\in\Lambda\setminus\{\lambda_1,\lambda_N\}$,
\bq
(\lambda-1-r)^2-4r&< & 0\eq
so that
\bq
\lambda&> & 1+r-2\sqrt{r}\ =\ (1-\sqrt{r})^2\eq
In particular, for $N>2$, we get
\bq
\lambda_2&> & (1-\sqrt{r})^2\eq
and as announced, for $N$ large,
\bq
\lambda_2-\lambda_1\ \sim \ \lambda_2\ \gg\ \lambda_1\eq
meaning that convergence to quasi-stationarity
happens at a much faster rate than absorption.\par
Theorem \ref{tvL2rev} shows that
for any $t\geq 0$,
\bq
\sup_{\mu_0\in\cP}\lVe \mu_t-\nu\rVe_{\mathrm{tv}}&\leq &
\sqrt{\frac{r^N}{r-1}}\lt(\frac{r}{r-1}\rt)^2(1+\cO(r^{-N}))\exp\lt(-(1-\sqrt{r}+\cO(r^{-N}))^2t\rt)
\eq
(where $\cO(r^{-N})$ is with respect to $N$, uniformly in $t\geq 0$).
It follows that for any fixed $s\geq 0$, if for $N$ large enough we consider the time
\bq
t&\df& \frac1{2(1-\sqrt{r})^2}\lt(\ln(r) N
+2s\rt)\eq
then
\bq
\sup_{\mu_0\in\cP}\lVe \mu_t-\nu\rVe_{\mathrm{tv}}&\leq &\frac{r^2}{(r-1)^{5/2}}\ (1+
\circ(1))\exp(-s)\eq
Notice  that the relaxation time to quasi-stationarity needs to be at least of order $N$, since it is already the order
of time required by the semi-group associated to $\wi L$ to get from $1$ to $N$, which supports a non-negligible
part of $\wi\eta$ (but starting from $N$, it can be shown that the relaxation time to quasi-stationarity is bounded
independently from $N$).

\subsection{A birth and death example with $\lambda_1\gg \lambda_2-\lambda_1$}\label{abadewlall}

The setting is as in the previous example, except that  $r<1$ in \eqref{L2} and \eqref{L2b}.\par
The beginning of Subsection \ref{abadewlll} is still valid: the reversible probability $\eta$ is 
given by \eqref{etar} and Lemma \ref{Lambda} is true, without the condition $N> (1+r)/(r-1)$, which is now void.
The difference starts with Proposition \ref{rgeq1}, which must be replaced by
\begin{pro}\label{rleq1}
For  $N\geq 4$,
we have
\bq
(1-\sqrt{r})^2+4\sqrt{r}\sin^2(
(1-r)/(2N+4))
\ \leq \ \lambda_1\ \leq \ (1-\sqrt{r})^2+4\sqrt{r}\sin^2(\pi/(2N))\\
 (1-\sqrt{r})^2+4\sqrt{r}\sin^2(\pi/(2N))\ \leq \ 
\lambda_2\ \leq \ (1-\sqrt{r})^2+4\sqrt{r}\sin^2(\pi/N)\eq
Furthermore, if $\varphi$ is an eigenvector associated to $\lambda_1$, we have
\bq
\frac{\varphi_\vee}{\varphi_\wedge}&\leq &r^{-(N-1)/2}\frac{1}{\sin((1-r)/(2N+4))}\eq
\end{pro}
\proof
First we show that none of the roots of $P_N$ is a real number.
Consider the function
\bq
f\st \RR_+\ni x\mapsto x^{N+1}-x^N+r^{1-N}x-r^{-N-1}\eq
According to the arguments of Lemma \ref{Lambda}, it is sufficient to show that $f$ only vanishes at $1/r$.
Its second derivative is given by
$f''(x)=(N+1)Nx^{N-1}-N(N-1)x^{N-2}$, for $x\geq 0$. Thus $f''$ is negative on  $(0,(N-1)/(N+1))$
and positive on $((N-1)/(N+1),+\iy)$. Furthermore, we compute that 
\bq
f'\lt(\frac{N-1}{N+1}\rt)&=&-\lt(1-\frac{2}{N+1}\rt)^{N-1}+r^{1-N}\eq
This quantity is positive for $N\geq 2$.
Thus $f$ is increasing on $(0,+\iy)$ and  can only vanish at $1/r$.\par
Since we know that the roots of $P_N$ are given by \eqref{rhopm} for $\lambda\in\Lambda\subset\RR$, we deduce that
\bq
\fo \lambda\in\Lambda,\qquad (\lambda-1-r)^2&<& 4r\eq
It follows that the modulus of $\rho_{\pm}$ in \eqref{rhopm} is given by $1/\sqrt{r}$, independently of $\lambda\in\Lambda$.
More precisely, there exists a set $\Theta\subset (0,\pi)$, such that the roots of $P_N$ are given by
\bq
\lt\{\frac1{\sqrt{r}}\exp(\pm i\theta)\st \theta\in \Theta\rt\}\eq
By using the mapping $\Psi$ of Lemma \ref{Lambda}, we get that  the spectrum of $L$ is 
\bq
\Lambda&=&\{l(\theta)\df 1+r-2\sqrt{r}\cos(\theta)\st \theta\in\Theta\}\eq
and that  corresponding eigenvectors are given by
\bqn{phitheta}
\fo x\in\lin 1,N\rin,\qquad \varphi_{\theta}(x)&=&r^{-x/2}\sin(\theta x)\eqn
for $\theta\in\Theta$ (note the slight modification of notation with respect to Lemma \ref{Lambda},
 indexing by  elements of $\Theta$ instead of $\Lambda$).
Ordering $\Theta$ into $0<\theta_1<\theta_2<\cdots < \theta_N<\pi$, it appears that
\bqn{llll}
\lambda_1=l(\theta_1)&\hbox{and}& \lambda_2=l(\theta_2)\eqn
From  Miclo \cite{MR2438701}, we deduce that $\varphi_{\theta_1}$ is non-decreasing 
and that $\varphi_{\theta_2}$ changes sign once (more generally 
$\varphi_{\theta_k}$ changes sign $k-1$ times, for $k\in\lin 1, N\rin$).
This remark and \eqref{phitheta} lead to the bounds
\bq
\theta_1\ \leq \ \pi/N&\hbox{and}&  \pi/N\ \leq \ \theta_2\ \leq \ 2\pi/N\eq
Taking into account that 
\bq
\fo \theta\in\Theta,\qquad 
l(\theta)&=&(1-\sqrt{r})^2+2\sqrt{r}(1-\cos(\theta))\\
&=&(1-\sqrt{r})^2+4\sin^2(\theta/2)\eq
and that sinus is positive and increasing on $(0,\pi/2)$,
we get that for $N\geq 4$,
\bq
\lambda_1\ \leq \ (1-\sqrt{r})^2+4\sqrt{r}\sin^2(\pi/(2N))\\
 (1-\sqrt{r})^2+4\sqrt{r}\sin^2(\pi/(2N))\ \leq \ 
\lambda_2\ \leq \ (1-\sqrt{r})^2+4\sqrt{r}\sin^2(\pi/N)\eq
To obtain a lower bound of the same kind for $\lambda_1$, recall that the elements $\theta\in\Theta$
satisfy the equation
\bqn{eqtheta}
\frac1{r^{N+1}}\exp(i2(N+1)\theta)-\frac1{r^N}\exp(i2N\theta)+\frac1{r^N}\exp(i2\theta)-\frac1{r^{N+1}}&=&0\eqn
and in particular $g(\theta)=0$, where the mapping $g$ is defined by
\bqn{eqtheta2}
\fo \theta\in\RR,\qquad g(\theta)&\df& \sin(2(N+1)\theta)-r\sin(2N\theta)+r\sin(2\theta)\eqn
One computes that $g'(0)=2N(1-r)+2+r$ and that
\bqn{eqtheta3}
\fo \theta\in\RR,\qquad \lve g''(\theta)\rve&\leq & 4(N+1)(N+2)\eqn
By consequence, the first zero of $g$ after 0 is larger than $(2N(1-r)+2+r)/(4(N+1)(N+2))$
and in particular
\bqn{theta1leq}
\nonumber\theta_1&\geq & \frac{2(1-r)(N+1)}{4(N+1)(N+2)}\\
&=&\frac{1-r}{2(N+2)}\eqn
leading to the announced lower bound on $\lambda_1$.\par
Furthermore, if 
$\varphi\df\varphi_{\theta_1}$,
we have
\bq
\frac{\varphi_\vee}{\varphi_\wedge}&=&\frac{ \varphi_{\theta_1}(N)}{ \varphi_{\theta_1}(1)}\\
&\leq & r^{-(N-1)/2}\frac{1}{\sin(\theta_1)}\\
&\leq&  r^{-(N-1)/2}\frac{1}{\sin((1-r)/(2N+4))}\eq
\wwtbp
\par
Working for fixed $r\in(0,1)$ in the asymptotic $N\ri\iy$, we deduce that
\bq
(1-\sqrt{r})^2+\sqrt{r}(1-r)^2\frac{1}{N^2}(1+\circ(1))
\ \leq \ \lambda_1\ \leq \ (1-\sqrt{r})^2+\frac{\pi^2}{N^2}(1+\circ(1))\\
 (1-\sqrt{r})^2+\frac{\pi^2}{N^2}(1+\circ(1))\ \leq \ 
\lambda_2\ \leq \ (1-\sqrt{r})^2+\frac{4\pi^2}{N^2}(1+\circ(1))\eq 
and
\bq
\frac{\varphi_\vee}{\varphi_\wedge}&\leq &\frac{2N}{(1-r)r^{(N-1)/2}}(1+\circ(1))\eq
In particular,  we get
\bq
\lambda_2-\lambda_1&\leq & \frac{4\pi^2-\sqrt{r}(1-r)^2}{N^2}(1+\circ(1))\\
\lambda_1&\sim & (1-\sqrt{r})^2\eq
and, as announced, for $N$ large,
\bq
\lambda_2-\lambda_1&\ll & \lambda_1\eq
meaning that  absorption
happens at a much faster rate than convergence to quasi-stationarity.\par
To exhibit an approximate quantitative estimate for the latter convergence, we need a lower bound 
on $\lambda_2-\lambda_1$.
\begin{lem}
For $N$ large enough, we have
\bq
\lambda_2-\lambda_1&\geq & \frac{(1-r)^2\sqrt{r}}{2N^2}(1+\circ(1))\eq
\end{lem}
\proof
With the notation of the proof of Proposition \ref{rleq1},
 begin by obtaining a lower bound on $\theta_2-\theta_1$.
Considering the function $g$ defined in \eqref{eqtheta2},  $\theta_2$ is larger
than the zero of $g$ following $\theta_1$.
We have
\bq
g'(\theta_1)&\df& 2(N+1)\cos(2(N+1)\theta_1)-r2N\cos(2N\theta_1)+2r\cos(2\theta_1)\eq
From \eqref{eqtheta}, we also obtain
\bq
\cos(2(N+1)\theta_1)-r\cos(2N\theta_1)+r\cos(\theta_1)-1&=&0\eq
so that
\bq
g'(\theta_1)&=&2N(1-r\cos(\theta_1))+2\cos(2(N+1)\theta_1)+2r\cos(2\theta_1)\\
&\geq & 2N(1-r)-2(1+r)\eq
Taking into account \eqref{eqtheta3}, we deduce that \bq
\theta_2-\theta_1&\geq & \frac{N(1-r)-(1+r)}{2(N+1)(N+2)}\\
&=& \frac{1-r}{2N}(1+\circ(1))
\eq
Next, we have
\bq
-\cos(\theta_2)&\geq &-\cos(\theta_1)+\min_{\theta\in[\theta_1,\theta_2]}\sin(\theta)(\theta_2-\theta_1)\\
&=& -\cos(\theta_1)+(1+\circ(1))\theta_1(\theta_2-\theta_1)\\
&\geq & -\cos(\theta_1)+\frac{1-r}{2(N+2)} \frac{1-r}{2N}(1+\circ(1))\eq
where we used \eqref{theta1leq}.
The announced bound is now a consequence of \eqref{llll}.
\wwtbp
\par
Putting together the previous estimates, with $\eta_\wedge \sim (1-r)r^{N-1}$, we get
\bq
\sup_{\mu_0\in\cP}\lVe \mu_t-\nu\rVe_{\mathrm{tv}}&\leq &\frac{4N^2}{(1-r)^{5/2}r^{3(N-1)/2}}(1+\circ(1))\exp\lt(
-\frac{(1-r)^2\sqrt{r}}{2N^2}(1+\circ(1))t
\rt)\eq
In particular, for any given $\epsilon >0$, if we consider 
\bq
t_N&\df& 4(1+\epsilon)\frac{N^2\ln(N)}{(1-r)^2\sqrt{r}}\eq
then
\bq
\lim_{N\ri\iy}\sup_{\mu_0\in\cP}\lVe \mu_{t_N}-\nu\rVe_{\mathrm{tv}}&=&0\eq

\subsection{A non-reversible example}\label{anre}

Let $N\in\NN$ be fixed. We consider $\bar S=S\sqcup \{\iy\}$, with $S=\ZZ_N$.
The generator $\bar L$ allows with rate 1 jumps of size 1 in $\ZZ_N$ and a jump at rate 1 from $0\in\ZZ_N$ to $\iy$, the absorbing point.
Namely, the generator $L$ is given by
\bq
\fo x,y\in \ZZ_N,\qquad L(x,y)&\df& \lt\{
\begin{array}{ll}
1&\hbox{, if $y=x+1$}\\
-1&\hbox{, if $y=x$}\\
0&\hbox{, otherwise}
\end{array}\rt.\eq
whose invariant probability measure $\eta$  is  the uniform distribution.
The potential $V$ takes the value 1 at 0 and 0 otherwise.
The spectral decomposition of the highly non-reversible operator $L-V$ is given by:
\begin{lem}\label{varphic}
Let $\cC$ be the set of (complex) solutions of the equation $X^N+X^{N-1}-1=0$.
Its cardinality is $N$ (i.e.\ all the solutions of the equation are distinct), the set of eigenvalues of
$L-V$ is $\{c-1\st c\in\cC\}$ and corresponding eigenvectors are given by the functions $\varphi_c$, for $c\in\cC$,
defined by
\bq
\fo x\in \lin 0, N-1\rin,\qquad \varphi_c(x)&\df& \lt\{
\begin{array}{ll}
1&\hbox{, if $x=0$}\\
c^{x-N}&\hbox{, otherwise}
\end{array}\rt.\eq
(where $\ZZ_N$ is naturally identified with $\lin 0, N-1\rin$).
\end{lem}
\proof
We begin by checking that all the roots of the polynomial  $X^N+X^{N-1}-1$ are simple.
Indeed, if $c\in\cC$ had multiplicity at least  two,  it would also
satisfy $Nc^{N-1}+(N-1)c^{N-2}=0$, namely $c=(1-N)/N$ (because 0 does not belong to $\cC$).
The equation $c^N+c^{N-1}=1$ could then be rewritten 
\bq
\frac1{N}\lt(\frac{1-N}{N}\rt)^{N-1}&=&1\eq
but this is impossible, because the absolute value of the l.h.s.\ is strictly less than 1.\par
Next we compute that for $c\in\cC$,
\bq
\fo x\in \ZZ_N,\qquad L[\varphi_c](x)&=& \lt\{
\begin{array}{ll}
c^{1-N}-1&\hbox{, if $x=0$}\\
(c-1)\varphi_c(x)&\hbox{, otherwise}
\end{array}\rt.\eq
Note that 
\bq
c^{1-N}-1&=&c^{1-N}-c+(c-1)\varphi_c(0)\\
&=&1+(c-1)\varphi_c(0)\\
&=&V(0)\varphi_c(0)+(c-1)\varphi_c(0)\eq
Thus it appears that on $\ZZ_N$,
\bq
(L-V)[\varphi_c]&=&(c-1)\varphi_c\eq
which is the wanted result, since we have exhibited exactly $N$ eigenvalues.\wwtbp
\par
Necessarily $\cC$ contains some real numbers, due to the Perron-Frobenius theorem which asserts that the smallest eigenvalue $\lambda_1$ of $V-L$  satisfies
\bq
\lambda_1&=&1-\max\{c\st c\in\cC\cap \RR\}\eq
By the strong irreversibility of $L$, the set $\cC\cap\RR$ is in fact very restricted, an observation which enables
 easy deduction of the asymptotic behavior of $\lambda_1$ for $N$ large:
\begin{lem}
If $N$ is odd, $\cC\cap\RR=\{1-\lambda_1\}$ and if $N$ is even, $\cC\cap\RR$ consists
 of two points. In both cases, 
 $\cC\cap\RR_+=\{1-\lambda_1\}$ and
 we have for $N$ large
 \bq
 \lambda_1&\sim& \frac{\ln(2)}{N}\eq
\end{lem}
\proof
Consider the function
\bq
g\st \RR\ni x&\mapsto& x^N+x^{N-1}-1\eq
The study of its variations leads to the two first announced results by differentiating it twice.
Indeed, if $N$ is odd, $g$ is increasing on $(-\iy,(1-N)/N)$, decreasing on $((1-N)/N,0)$
and increasing on $(0,+\iy)$. As  was already seen in the proof
of the previous lemma, $g((1-N)/N)<0$, so that $g$ admits a unique real root  contained in $(0,+\iy)$.
For $N$ even, 
$g$ is decreasing on $(-\iy,(1-N)/N)$
and increasing on $((1-N)/N,+\iy)$. Since $g((1-N)/N)<0$ and $\lim_{\pm\iy}g=+\iy$, $g$ admits two real roots,
the largest one being the unique one belonging to $(0,+\iy)$, since $g(0)=-1$.
\par
Let $y>0$ be given and for $N>y$ consider $x_N=1-y/N$. It appears that
\bq
\lim_{N\ri\iy}g(x_N)&=&2\exp(-y)-1\eq
It follows that the unique root $c_N$ of $g$ in $(0,+\iy)$ satisfies for $N$ large
\bq
c_N-1&\sim&- \frac{\ln(2)}{N}\eq
which amounts to the last announced result.
\wwtbp
\par
Let $\varphi=\varphi_{1-\lambda_1}$, with the notation of Lemma \ref{varphic}, be an eigenvector associated to $\lambda_1$.
We have, for $N$ large
\bq
\frac{\varphi_\vee}{\varphi_\wedge}
&=&\frac{\varphi(1)}{\varphi(0)}\\
&=&(1-\lambda_1)^{1-N}\\
&\sim& \exp(\ln(2))\ =\ 2\eq
\par
In addition, note that $L^*$, the dual operator of $L$ in $\LL^2(\eta)$, is given by 
\bq
\fo x,y\in \ZZ_N,\qquad L^*(x,y)&\df& \lt\{
\begin{array}{ll}
1&\hbox{, if $y=x-1$}\\
-1&\hbox{, if $y=x$}\\
0&\hbox{, otherwise}
\end{array}\rt.\eq
It corresponds to the conjugation of $L$ with the involutive transformation of $\ZZ_N$
given by $\iota\st\ZZ_N\ni x\mapsto -x$ (or $\lin 1,N-1\rin\ni x\mapsto N-x$ and $\iota(0)=0$). It follows that the function $\varphi^*$ considered
in the introduction is proportional to $\varphi\circ \iota$, so that the mapping $\varphi\varphi^*$ is constant.
In particular the probability $\wi\eta$ defined in \eqref{wieta} is equal to $\eta$, the uniform distribution on $\ZZ_N$.
Furthermore, we compute that
the generator $\wi L$ defined in \eqref{wiL} is given by
\bq
\fo x,y\in \ZZ_N,\qquad \wi L(x,y)&\df& \lt\{
\begin{array}{ll}
(1-\lambda_1)&\hbox{, if $x\not= 0$ and $y=x+1$}\\
-(1-\lambda_1)&\hbox{, if $x\not= 0$ and $y=x$}\\
(1-\lambda_1)^{1-N}&\hbox{, if $x=0$ and $y=1$}\\
-(1-\lambda_1)^{1-N}&\hbox{, if $x=0$ and $y=0$}\\
0&\hbox{, otherwise}
\end{array}\rt.\eq
Its additive symmetrization $\wi L^\diamond$  in $\LL^2(\eta)$ gives the rate $(1-\lambda_1)/2$ to any oriented edge $(x,x+1)$ or $(x+1,x)$ of $\ZZ_N$,
except to the edges $(0,1)$ and $(1,0)$, which have the rate $(1-\lambda_1)^{1-N}/2$.
By comparison with the usual continuous-time random walk on $\ZZ_N$, we deduce that the spectral gap $\wi \lambda$ of $\wi L^\diamond$
satisfies
\bq
(1-\cos(2\pi/N))(1-\lambda_1)\ \leq \ \wi\lambda\ \leq \ (1-\cos(2\pi/N))(1-\lambda_1)^{1-N}\eq
namely, asymptotically for $N$ large,
\bq
\frac{2\pi^2}{N^2}(1+\circ(1))\ \leq \  \wi\lambda\ \leq \ \frac{4\pi^2}{N^2}(1+\circ(1))\eq
\par
Relying on \eqref{will}, 
we would have obtained
\bq
\wi \lambda&\geq & \frac{1+\circ(1)}{4}\lambda\eq
where $\lambda$ is the spectral gap of the additive symmetrization of $L$ in $\LL^2(\eta)$, 
which is the usual continuous-time random walk on $\ZZ_N$, so that $\lambda\sim 2\pi^2/N^2$.
Thus it only leads to a slight deterioration on the estimate of $\wi\lambda$  obtained by working directly with
\eqref{Poincare}.\par
For large $N$, Theorem \ref{tvL2} leads to
\bq
\sup_{\mu_0\in\cP}\lVe \mu_t-\nu\rVe_{\mathrm{tv}}&\leq &
2\sqrt{N}(1+\circ(1))\exp\lt(\frac{2\pi^2}{N^2}(1+\circ(1))t\rt)\eq
In particular, for any given $\epsilon >0$, if we consider 
\bq
t_N&\df& (1+\epsilon)\frac{N^2\ln(N)}{4\pi^2}\eq
then
\bq
\lim_{N\ri\iy}\sup_{\mu_0\in\cP}\lVe \mu_{t_N}-\nu\rVe_{\mathrm{tv}}&=&0\eq

\subsection{A product example}\label{ape}

Let us first come back to the general setting of the introduction (which is then tensorized).
Let $d\in\NN$, be given. On $ S^d$,  consider the Markovian generator
\bq
 L^{(d)}&\df& \frac1d\sum_{k\in\lin 1,d\rin}  L_k\eq
where $ L_k$ acts like $ L$ on the $k$-th coordinate of $S^d$.
Define furthermore  the potential  $V^{(d)}$ by
\bq
\fo x\df(x_1, ..., x_d)\in S^d,\qquad V^{(d)}&\df&  \frac1d\sum_{k\in\lin 1,d\rin} V(x_k)\eq
Note that the associated $\bar L^{(d)}$ is not of the form  $(1/d)\sum_{k\in\lin 1,d\rin}  \bar L_k$, because 
the underlying state space would be $(\bar S)^d$
and not $S^d\sqcup\{\iy\}$ as it should be. One recovers the subMarkovian generator $L^{(d)}-V^{(d)}$
by modifying $(1/d)\sum_{k\in\lin 1,d\rin}  \bar L_k$ so that  all the points of $\{x\df(x_1, ..., x_d)\in (\bar S)^d\st \ex k\in\lin 1, d\rin \hbox{ with } x_k=\iy\}$
become absorbing.
\par
The invariant measure $\eta^{(d)}$ associated to $ L^{(d)}$ is $\eta^{\otimes d}$
and we have
\bq
L^{(d)}-V^{(d)}&=& \frac1d\sum_{k\in\lin 1,d\rin} ( L-V)_k\eq
It appears in particular that the  first eigenvalue of $V^{(d)}-L^{(d)}$ is $\lambda_1$,
the same as that of $V-L$ and the associated quasi-stationary distribution (respectively
first eigenfunction) is $\nu^{\otimes d}$ (resp.\ $\varphi^{\otimes d}$). 
It follows that $\wi L^{(d)}$, the Doob transform of $L^{(d)}-V^{(d)}$ by $\varphi^{\otimes d}$,
satisfies
\bq\wi L^{(d)}&=& \frac1d\sum_{k\in\lin 1,d\rin} \wi L_k\eq
and that its invariant probability $\wi\eta^{(d)}$ is  $\wi\eta^{\otimes d}$.
 In a similar way, we have that
\bq L^{*(d)}&=& \frac1d\sum_{k\in\lin 1,d\rin} L^*_k\eq
and the  first eigenvector of  $- L^{*(d)}$ is $(\varphi^*)^{\otimes d}$.
Finally  $\wi L^{\diamond(d)}$, the additive symmetrization of $\wi L^{(d)}$ in $\LL^2(\wi\eta^{\otimes d})$, is equal
to $(1/d)\sum_{k\in\lin 1,d\rin} \wi L^{\diamond}_k$, so that its spectral gap $\wi \lambda$ (respectively its logarithmic Sobolev constant $\wi\alpha$) is equal 
to that of $ \wi L^{\diamond}$ (for such tensorization properties, see for instance the book \cite{MR2002g:46132} of Ané et al.).
\par
With obvious notation, Theorem \ref{tvL2}
then leads to the fact that for any $t\geq 0$, we have
\bq
\sup_{\mu_0^{(d)}\in\cP^{(d)}}\lVe \mu_t^{(d)}-\nu^{\otimes d}\rVe_{\mathrm{tv}}&\leq &
\lt(\sqrt{\frac{\eta[\varphi\varphi^*]}{(\varphi\varphi^*\eta)_{\wedge}}}
\frac{\varphi_\vee}{\varphi_\wedge}\rt)^d\exp(-\wi\lambda t)\eq
Under the reversibility condition of Theorem 
\ref{tvL2rev},
we get that
 for any $t\geq 0$, 
\bqn{L22}
\sup_{\mu_0^{(d)}\in\cP^{(d)}}\lVe \mu_t^{(d)}-\nu^{\otimes d}\rVe_{\mathrm{tv}}&\leq &
\lt( \sqrt{\frac{1}{\eta_{\wedge}}}
\lt(\frac{\varphi_\vee}{\varphi_\wedge}\rt)^2\rt)^d\exp(-(\lambda_2-\lambda_1) t)
\eqn
The bound \eqref{ent} can be rewritten in the form
\bq
\sup_{\mu_0^{(d)}\in\cP^{(d)}}\lVe \mu_t^{(d)}-\nu^{\otimes d}\rVe_{\mathrm{tv}}&\leq &
\sqrt{2d\ln\lt(\frac{\eta[\varphi\varphi^*]}{(\varphi\varphi^*\eta)_{\wedge}}\rt)
\lt(\frac{\varphi_\vee}{\varphi_\wedge}\rt)^d}\exp(-(\wi\alpha/2) t)\eq
or under the reversibility condition,
\bqn{ent2}
\sup_{\mu_0^{(d)}\in\cP^{(d)}}\lVe \mu_t^{(d)}-\nu^{\otimes d}\rVe_{\mathrm{tv}}&\leq &
\sqrt{2d\ln\lt(\frac{1}{\eta_{\wedge}}\frac{\varphi_\vee}{\varphi_\wedge}\rt)
\lt(\frac{\varphi_\vee}{\varphi_\wedge}\rt)^d}\exp(-(\wi\alpha/2) t)
\eqn
\par
It is easy to construct an example showing that \eqref{ent2} can lead to a better estimate than \eqref{L22}.
Take 
\bq S\df\{1,2\},\quad L\df\lt(\begin{array}{cc} -1&1\\
1&-1\end{array}\rt),\quad V\df\lt(\begin{array}{c} 1\\
1\end{array}\rt)\eq
for which $\varphi\equiv 1$, $\eta=(1/2,1/2)$, $\lambda_2-\lambda_1=2$ and $\wi\alpha =1$
(recall the convention after \eqref{alla}, which is an equality in this two-points case, see Diaconis and Saloff-Coste \cite{MR1410112}).
The r.h.s.\ of \eqref{L22} and \eqref{ent2} are respectively $2^{d/2}\exp(-2t)$ and $\sqrt{2d\ln(2)}\exp(-t/2)$.
The first bound leads to a mixing time (the first time $t>0$ the quantity
$\sup_{\mu_0^{(d)}\in\cP^{(d)}}\lVe \mu_t^{(d)}-\nu^{\otimes d}\rVe_{\mathrm{tv}} $ goes below a fixed level such as $1/2$)
of order $d$, while the second bound rather gives order $\ln(d)$.
\par
For a little less artificial example, one can come back to Subsection \ref{abadewlall}, with $N=2$ and $r>0$ very small.
Indeed, one computes that
\bq
\lambda_2-\lambda_1&=&2\sqrt{r(1+r)}\\
\frac{\varphi_\vee}{\varphi_\wedge}&=&\sqrt{\frac{1-r}{r}}\\
\eta_\wedge&=&\frac{r}{1+r}\\
\wi\eta_\wedge&=&\frac{1-r}{2}\eq
It follows from \eqref{alla} that for $0<r\ll 1$,
\bq
\wi\alpha& \geq & \frac{r}{\ln((1+r)/(1-r))}(\lambda_2-\lambda_1)\\
&\sim& \frac{\lambda_2-\lambda_1}{2}\eq
For  $r>0$ small,
we get from \eqref{L22} and \eqref{ent2} that the leading term in $d\in\NN$ in  the deduced upper bounds on the mixing time are respectively $3d/(4\sqrt{r})\ln(1/r)$
and $d/(2\sqrt{r})\ln(1/r)$, showing thus a little advantage for the estimate coming from \eqref{ent2}.

\section{Some discrete time models}\label{sdtm}
Of course the theory can be developed in discrete time as well. We briefly carry this out here and treat some higher dimensional examples where all the spectral information is available.
Let $\bar S\df S\sqcup\{\iy\}$ be the extended state space with $\iy$ the absorbing state. Denote by $N$ the cardinality of $S$.
The transition matrix can be written
\bq
\lt(\begin{array} {cccc}
1& 0&\quad\cdots\quad& 0\\
a_1& & & \\
 & & & \\
 \vdots& &Q&\\
  & & & \\
a_N& & & \\
\end{array}
\rt)\eq
with $Q$ an $N\times N$ matrix, here assumed to be irreducible. 
Let $\psi$ and $\varphi$ be positive left and right eigenvectors of $Q$ with eigenvalue $\beta >0$ of largest size.
Set 
\bq
\fo x,y\in S,\qquad K(x,y)&\df& Q(x,y)\frac{\varphi(x)}{\beta\varphi(y)}\eq
This is a Markov transition matrix on $S$ with stationary distribution $\pi$ given by
\bq
\fo x\in S,\qquad \pi(x)&\df& \frac{\varphi(x)\psi(x)}{\sum_{y\in S} \varphi(y)\psi(y)}\eq
It has the probabilistic interpretation of the transition probabilities for the original chain conditioned on non-absorption (for all time).
The quasi-stationary distribution is given by
\bq
\fo x\in S,\qquad \nu(x)&\df& \frac{\psi(x)}{\sum_{y\in S} \psi(y)}\eq
observe that the ratio $r\df \varphi_\vee/\varphi_\wedge$ allows the bounds 
\bq
\fo x\in S,\qquad r^{-1} \nu(x)\ \leq \ \pi(x)\ \leq \ r\nu(x)\eq
\par
If $Q$ above is diagonalizable, with right 
eigenfunctions $(f_i)_{i\in\lin N\rin}$ and left  
eigenfunctions $(g_i)_{i\in\lin N\rin}$
for  eigenvalues $(\beta_i)_{i\in\lin N\rin}$, 
normalized so that $\sum_{x\in S}g_i(x)f_j(x)=\delta_{i,j}$ for any $i,j\in\lin N\rin$,
then 
\bq
\fo l\in\ZZ_+,\,\fo x,y\in S,\qquad
Q^l(x,y)&=&\sum_{k=1}^N\beta_k^lf_k(x)g_k(y)\eq
thus 
\bq
P[X_l=y\vert X_0=x, T>l]&=&\frac{Q^l(x,y)}{Q^l(x)}\eq
with $Q^l(x)=\sum_{y\in S} Q^l(x,y)$ and where $(X_n)_{n\in\ZZ_+}$ is the underlying absorbing Markov chain and $T$ is its absorbing time.
\par
Explicit diagonalizations are available surprisingly often. For example, for a birth and death chain on $\{0, 1, ..., 2N\}$, symmetric with respect to $N$, take the starting point to be zero 
and the absorbing point to be $N$. If $(\varphi_i)_{i\in\lin 0, 2N\rin}$ are the right eigenvectors of the original chain, often available as orthogonal polynomials,
$\varphi_1$, $\varphi_3$, ..., $\varphi_{2N-1}$ all vanish at $N$ and so  restrict to the needed $(f_i)_{i\in\lin N\rin}$.
Because birth and death chains are reversible, these determine the family $(g_i)_{i\in\lin N\rin}$ and the ingredients for analysis are available. The Ehrenfest urn and the example at the end of this section are two cases where we have carried this approach out to get sharp answers (matching upper and lower bounds for convergence to quasi-stationarity). It is only fair to report that the analysis involved can require substantial effort.
\subsection{Example of rock breaking}
In this example the matrix $Q$ is not irreducible, nevertheless the above results can be applied, because the function $\varphi$ is (strictly) positive.
To justify this observation, replace for $\epsilon\in(0,1)$, $Q$ by $(1-\epsilon)Q+\epsilon J$, where $J$ has all its entries equal to $1/N$, apply the previous results and let $\epsilon$ go to zero.\par
Let $n\in\NN$ be given and $\bar S\df\cP(n)$, the set of all partitions of $n$. Thus if $n=4$,
$\bar S=\{ 4, 31, 22, 211, 1111\}$. An absorbing Markov chain on $\bar S$, modeled on a rock breaking Markov chain studied by Kolmogorov, is developed in 
Diaconis, Pang and Ram \cite{2012arXiv1206.3620D}. Briefly, if $\lambda=(\lambda_1, \lambda_2, ...,\lambda_l)$, with $\lambda_1\geq \lambda_2\geq \cdots \geq  \lambda_l>0$,
$\lambda_1+\cdots +\lambda_l=n$, the chain proceeds from $\lambda$ by independently choosing, for $i\in\lin 1, l\rin$, binomial variables $\lambda_i^{(1)}$ of parameters $(\lambda_i, 1/2)$,
so that we can write $\lambda_i\fd \lambda_i^{(1)}+\lambda_i^{(2)}$. Next, after discarding any zeros and reordering the  $\lambda_i^{(1)},\lambda_i^{(2)}$, for $i\in\lin 1, l\rin$, we get the 
new position of the chain. It is absorbing at $(1^n)$. The natural starting place is $(n)$.
\par
In \cite{2012arXiv1206.3620D}, the eigenvalues are shown to be $1, 1/2, 1/4, ..., 1/2^n$, with $1/2^{n-l}$ having multiplicity $p(n,l)$, the number of partitions of $n$ into $l$ parts.
In particular the second eigenvalue is $1/2$, with multiplicity 1. The eigenvectors are given explicitly and these restrict to give explicit left and right eigenbases of $Q$.
With notation as above, for $\beta=1/2$, for all $\lambda=(\lambda_1, \lambda_2, ..., \lambda_l)\in\cP(n)$,
\bq
\varphi(\lambda)&=&\sum_{i\in\lin 1, l\rin}\binom{\lambda_i}{2}\\
\psi(\lambda)&=&\lt\{\begin{array}{ll}
1&\hbox{, if }\lambda=(1^{n-2},2)\\
0&\hbox{, otherwise}\end{array}\rt.
\eq
Thus $ \varphi_\vee/\varphi_\wedge=\binom{n}{2}/1=\binom{n}{2}$.
When $n=4$, the original transition matrix is
\bq
\begin{array}{cc}
 &\begin{array}{ccccc}1^4 &1^22& 2^2  & 13& 4\end{array}\\
 \begin{array}{c}1^4 \\1^22 \\ 2^2  \\13 \\ 4\end{array}&
\lt( \begin{array}{ccccc}1 &0& 0 &0& 0\\
1/2&1/2&0&0&0\\
1/4&1/2&1/4&0&0\\
0&3/4&0&1/4&0\\
0&0&3/8&1/2&1/8
\end{array}\rt)
\end{array}
\eq
The left (right) eigenvectors are given as the rows (columns) of the two arrays
\bq
\lt( \begin{array}{ccccc}1 &0& 0 &0& 0\\
1&1&0&0&0\\
1&2&1&0&0\\
1&3&0&1&0\\
1&6&3&4&0
\end{array}\rt)
\qquad\qquad
\lt( \begin{array}{ccccc}1 &0& 0 &0& 0\\
-1&1&0&0&0\\
1&-2&1&0&0\\
2&-3&0&1&0\\
-6&12&-3&-4&1
\end{array}\rt)\eq
So $\psi=(1,0,0,0)$, $\varphi=(1,2,3,6)^{\mathrm{t}}$. The adjusted transition matrix $K$ is given by
\bq
\begin{array}{cc}
 &\begin{array}{cccc} 1^22& 2^2  & 13& 4\end{array}\\
 \begin{array}{c} 1^22 \\ 2^2  \\13 \\ 4\end{array}&
\lt( \begin{array}{cccc}1 &0& 0 &0\\
1/2&1/2&0&0\\
1/2&0&1/2&0\\
0&1/4&1/2&1/4\\
\end{array}\rt)
\end{array}
\eq
The reader may check that discarding the top row and first column of the eigenvector arrays gives the eigenvectors of $K$.
\par
In this example the quasi-stationary distribution $\nu$ is the stationary distribution $\pi$ of $K$, both are the Dirac mass  at $1^{n-2}2$. The chain $K$ is itself absorbing.
This rock breaking chain is a special case of a host of explicitly diagonalizable Markov chains derived from Hopf algebras \cite{2012arXiv1206.3620D}. Some other algebraic constructions leading to explicit quasi-stationary calculations may be found in Defosseux \cite{2013arXiv1307.3830D} (fusion coefficients and random walks in alcoves of affine Lie algebras). Symmetric function theory, in various deformations (Sekiguchi-Debiard operators) leads to further explicit diagonalizations in the work of Jiang \cite{2012arXiv1204.1671J}. Turning either of these last sets of examples into sharp bounds seems like a fascinating research project.
\subsection{Geometric theory}
The basic path arguments of Holley and Stroock \cite{MR90g:60091}, Jerrum and Sinclair \cite{MR1025467} and Diaconis and Stroock \cite{MR1097463} can be applied to absorbing chains.
This was done in a sophisticated context in Miclo \cite{miclo:hal-00777146}. The following paragraph develops a simple version in the discrete context.
Let $S$ be a finite set, $\iy$ an absorbing point and $K$ a Markov chain on $\bar S\df S\sqcup\{\iy\}$. We suppose as above that the chain is absorbing with probability one and that the chain restricted to $S$ is connected. Suppose that $q$ is a probability on $S$ and consider $\LL^2(q)$ endowed with its usual inner product $\lan f,g\ran_q\df\sum_{x\in S} f(x)g(x)q(x)$, for $f,g\in \LL^2(q)$.
Suppose too that $q(x)K(x,y)=q(y)K(y,x)$ for $x,y\in S$. When needed, define $q(\iy)=0$ and the functions from $\LL^2(q)$
are extended on $\bar S$ by making them vanish at $\iy$. Let $\beta_1$ be the largest eigenvalue of $K$ restricted to $S$. The minimax characterization gives
\begin{lem}\label{Plem1}
If the Poincaré inequality $\lVe f\rVe^2_q\leq A\lan (I-K)f,f\ran_q$ holds for all $f\in \LL^2(q)$, then $\beta_1\leq 1-1/A$.
\end{lem}
\begin{rem} Of course, $\lambda_1$ of Section \ref{intro} satisfies $\lambda_1=1-\beta_1$.\end{rem}
\par
Define a Dirichlet form $\cE$ on $\LL^2(q)$, by
\bq
\fo f\in\LL^2(q),\qquad \cE(f,f)&\df& \frac12\sum_{x,y\in \bar S}(f(y)-f(x))^2\, q(x)K(x,y)\eq
\begin{lem}
For $f\in\LL^2(q)$, we have
\bq
\cE(f,f)&=& \lan (I-K)f,f\ran_q-\frac12\sum_{x\in S} f^2(x) q(x)K(x,\iy)\\
&\leq &\lan (I-K)f,f\ran_q\eq
\end{lem}
\proof 
This is simple by directly computing both sides of the equality, separating the cases where $x,y\in S$
and the cases where at least one of them is equal to $\iy$.\wwtbp
\par
To bring in geometry, for $x\in S$, let $\gamma_x$ be a path starting at $x$ and ending at $\iy$ with steps possible with respect to $K$. If there are many absorbing points,
$\gamma_x$ may connect $x$ to any of them. Thus $\gamma_x=(x_0=x, x_1, ..., x_l=\iy)$ with $K(x_i,x_{i+1})>0$ for $0\leq i\leq l-1$.
Let the length $l$ of the path be denoted $\lve \gamma_x\rve$.
\begin{pro}\label{geom}
With the notation as above, $A$ in Lemma \ref{Plem1} may be taken as
\bq A&=& \max_{x\in S, \, y\in \bar S}\frac2{q(x)K(x,y)}\sum_{z\in S\st (x,y)\in \gamma_z}
\lve \gamma_z\rve q(z)\eq
\end{pro}
\proof
Let $x\in S$ be given and write $\gamma_x=(x_0, x_1, ..., x_l)$. The idea is to expand
\bq
f^2(x)&=&(f(x_0)-f(x_1))+(f(x_1)-f(x_2) )+\cdots+ (f(x_{l-1})-f(x_l)))^2\\
&\leq & \lve \gamma_x\rve\sum_{i=0}^{l-1}(f(x_i)-f(x_{i+1}))^2\eq
Thus for $f\in \LL^2(q)$,
\bq
\sum_{x\in S} f^2(x)q(x)&=&\sum_{x\in S}
\lt(\sum_{e\in \gamma_x}f(e^-)-f(e^+)\rt)^2 q(x)\\
&\leq & \sum_{x\in S}\lve \gamma_x\rve\,q(x)\sum_{e\in\gamma_x}(f(e^-)-f(e^+))^2
\\
&=& \sum_{x\in S,\, y\in\bar S}(f(x)-f(y))^2\frac{q(x)K(x,y)}{q(x)K(x,y)}\sum_{z\in S\st (x,y)\in \gamma_z}
\lve \gamma_z\rve q(z)\\
&\leq & A\cE(f,f)\eq
\wwtbp
\begin{rem}
Path technology has evolved: with many choices of paths, one may choose randomly, see Diaconis and Saloff-Coste
\cite{MR1385408}, weights may be used in the Cauchy-Schwarz bound, as in Diaconis and Saloff-Coste
\cite{MR1649805}. This can be important when the stationary distributions varies a lot. Paths may be used locally, see Diaconis and Saloff-Coste
\cite{MR1385408}. Any such variation is easy to adapt in the above proposition.
\end{rem}
\subsection{Other examples in discrete time}
The following calculations are classical. The neat form presented here is borrowed from the thesis work of Hua Zhou \cite{Zhou_thesis}
and provides an alternative approach to Examples \ref{abadewlsll}, \ref{abadewlll} and \ref{abadewlall}.
Consider a birth and death chain on $\lin 0, N\rin$ with transition matrix
\bq
\lt(
\begin{array}{ccccc}
r&1-r& &&\\
p&0&q &&\\
&&\ddots&&\\
&&p&0&q\\
&&&1-s&s
\end{array}\rt)\eq
with $p\in(0,1)$, $q=1-p$, $r,s\in [0,1]$.
\begin{pro}
The eigenvalues and right eigenfunctions are of form
\bq
\beta \df2\sqrt{pq}\cos(\theta),\qquad \fo x\in\lin 0,N\rin,\quad \varphi(x) \df \lt(\frac{p}{q}\rt)^{x/2}\cos(\theta x+c)\eq
where $\theta$ and $c$ are determined by the boundary values:
\bq
r\varphi(0)+(1-r)\varphi(1)&=&\beta \varphi(0)\\
(1-s)\varphi(N-1)+s\varphi(N)&=&\beta \varphi(N)\eq
\end{pro}
\proof
This follows from the trigonometric identity 
\bq
\fo \alpha,\beta \in\RR,\qquad \cos(\alpha)+\cos(\beta)&=&2\cos\lt(\frac{\alpha+\beta}{2}\rt)\cos\lt(\frac{\alpha-\beta}{2}\rt)\eq
Since for any $\theta, c\in\RR$,
\bq
q \lt(\frac{p}{q}\rt)^{(x+1)/2}\cos(\theta (x+1)+c)+p \lt(\frac{p}{q}\rt)^{(x-1)/2}\cos(\theta (x-1)+c)&=&2\sqrt{pq}
 \lt(\frac{p}{q}\rt)^{x/2}\cos(\theta x+c)\eq
 Zhou \cite{Zhou_thesis} has shown that the above boundary conditions lead indeed to $N+1$ eigenvalues.
 \wwtbp
 As an example, take $p=q=1/2=r$, $s=1$. This gives the simple random walk on $\lin 0,N\rin$ absorbing at $N$
 with holding at 0.
 The above proposition gives the equations
 \bq
 \cos(c)+\cos(\theta +c)\ =\ 2\cos(\theta)\cos(c),\quad
 \cos(N\theta +c)\ =\ 0\ \hbox{ or }\  \theta\ =\ 0\eq
 These have solutions $c=\theta/2$, $\theta=j\pi/(2N-1)$, for $j=0,1, 3, ..., 2N-1$.
 It follows that the chain has eigenvalues $\beta_j\df \cos(j\pi/(2N-1))$, for $j=0,1, 3, ..., 2N-1$
 with right eigenfunctions $\varphi_j$ given by
 \bq
 \fo x\in \lin 0, N\rin,\qquad
\varphi_j(x)&\df& \cos\lt(\frac{(2x+1)j\pi}{2(2N+1)}\rt)\eq
The left eigenfunctions are $\psi_0(x)=\delta_N(x)$ and for $j=1, 3, ..., 2N-1$,
\bq
\fo x\in\lin 0, N\rin,\qquad
\psi_j(x)&\df& \lt\{
\begin{array}{ll}
\varphi_j(x)&\hbox{, if $x\in\lin 0, N-1\rin$}\\
\frac{(-1)^{(j+1)/2}}{2}\cot\lt(\frac{j\pi}{2(2N-1)}\rt)
&\hbox{, if $x=N$}\end{array}\rt.\eq
In particular, the quasi-stationary distribution has probability density given by
\bq
\fo x\in\lin 0, N-1\rin,\qquad \nu(x)&\df&\frac{\psi_1(x)}{\sum_{y\in\lin 0,N-1\rin}\psi_1(y)}
\\
&=& 2\tan\lt(\frac{\pi}{2(2N-1)}\rt)\cos\lt(\frac{(2x+1)\pi}{2(2N-1)}\rt)
\eq
\par
Consider the geometric bound from Proposition \ref{geom}. From the discussion above,
\bq
\beta_1\ =\ \cos\lt(\frac{\pi}{2N-1}\rt)\ =\ 1-\frac{\pi^2}{2(2N-1)^2}+\cO\lt(\frac1{N^4}\rt)\eq
The reversing probability $q$ on $\lin 0,N-1\rin$ is the uniform distribution.
There is a unique choice of (not self-intersecting) paths from $x$ to $N$. The quantity $A$ is obviously maximized at the edge $(N-1,N)$.
Then, it is
\bq
A\ =\ \frac{4N}{N}\sum_{x\in\lin 0, N-1\rin}N-x\ =\ 2N(N+1)\eq
This gives $\beta_1\leq 1-\frac{1}{4N(N+1)}$ which compares reasonably with the correct answer.\par\sm
In this problem, $\varphi_\vee/\varphi_\wedge$ is of order $N^2$ and our bounds show that order $N^2\ln(N^2)$ steps
suffice for convergence to quasi-stationarity. Using all of the spectrum, classical analysis shows that order $N^2$ steps are necessary and sufficient.
Zhou \cite{Zhou_thesis} gives similar exact formulae for reflecting and absorbing boundaries at zero. He also derives the exact spectral data for some absorbing birth and
death chains from biology (Morans model with various types of mutation).
 \bigskip
 \par
 \textbf{\large Acknowledgments:}\par\sm\noindent 
 We are grateful to Bertrand Cloez for discussions on the subject, as well as for pointing out some related references.
 Thanks also to Hua Zhou for help with Section \ref{sdtm} and to the ANR STAB for its support.

\vskip2cm
\hskip70mm\box5

\end{document}